# MAXIMUM $LQ$-LIKELIHOOD ESTIMATION

### By Davide Ferrari and Yuhong Yang[1]

*Università di Modena e Reggio Emilia and University of Minnesota*

In this paper, the maximum L$q$-likelihood estimator (ML$q$E), a new parameter estimator based on nonextensive entropy [*Kibernetika* **3** (1967) 30–35] is introduced. The properties of the ML$q$E are studied via asymptotic analysis and computer simulations. The behavior of the ML$q$E is characterized by the degree of distortion $q$ applied to the assumed model. When $q$ is properly chosen for small and moderate sample sizes, the ML$q$E can successfully trade bias for precision, resulting in a substantial reduction of the mean squared error. When the sample size is large and $q$ tends to 1, a necessary and sufficient condition to ensure a proper asymptotic normality and efficiency of ML$q$E is established.

**1. Introduction.** One of the major contributions to scientific thought of the last century is information theory founded by Claude Shannon in the late 1940s. Its triumph is highlighted by countless applications in various scientific domains including statistics. The central quantity in information theory is a measure of the "amount of uncertainty" inherent in a probability distribution (usually called Shannon's entropy). Provided a probability density function $p(x)$ for a random variable $X$, Shannon's entropy is defined as $\mathcal{H}(X) = -E[\log p(X)]$. The quantity $-\log p(x)$ is interpreted as the information content of the outcome $x$, and $\mathcal{H}(X)$ represents the average uncertainty removed after the actual outcome of $X$ is revealed. The connection between logarithmic (or additive) entropies and inference has been copiously studied (see, e.g., Cover and Thomas [9]). Akaike [3] introduced a principle of statistical model building based on minimization of entropy. In a parametric setting, he pointed out that the usual inferential task of maximizing the log-likelihood function can be equivalently regarded as minimization of the empirical version of Shannon's entropy, $-\sum_{i=1}^{n} \log p(X_i)$. Rissanen proposed

Received November 2007; revised January 2009.

[1]Supported in part by NSF Grant DMS-07-06850.

the well-known minimum description length criterion for model comparison (see, e.g., Barron, Rissanen and Yu [5]).

Since the introduction of Shannon's entropy, other and more general measures of information have been developed. Rényi [27] and Aczél and Daróczy [2] in the mid-1960s and 1970s proposed generalized notions of information (usually referred to as Rényi entropies) by keeping the additivity of independent information, but using a more general definition of mean. In a different direction, Havrda and Charvát [16] proposed *nonextensive* entropies, sometimes referred to as *q*-order entropies, where the usual definition of mean is maintained while the logarithm is replaced by the more general function $L_q(u) = (u^{1-q} - 1)/(1-q)$ for $q > 0$. In particular, when $q \to 1$, $L_q(u) \to \log(u)$, recovering the usual Shannon's entropy.

In recent years, *q*-order entropies have been of considerable interest in different domains of application. Tsallis and colleagues have successfully exploited them in physics (see, e.g., [29] and [30]). In thermodynamics, the *q*-entropy functional is usually minimized subject to some properly chosen constraints, according to the formalism proposed by Jaynes [19] and [20]. There is a large literature on analyzing various loss functions as the convex dual of entropy minimization, subject to constraints. From this standpoint, the classical maximum entropy estimation and maximum likelihood are seen as convex duals of each other (see, e.g., Altun and Smola [4]). Since Tsallis' seminal paper [29], *q*-order entropy has encountered an increasing wave of success and Tsallis' nonextensive thermodynamics, based on such information measure, is nowadays considered the most viable candidate for generalizing the ideas of the famous Boltzmann–Gibbs theory. More recently, a number of applications based on the *q*-entropy have appeared in other disciplines such as finance, biomedical sciences, environmental sciences and linguistics [14].

Despite the broad success, so far little effort has been made to address the inferential implications of using nonextensive entropies from a statistical perspective. In this paper, we study a new class of parametric estimators based on the *q*-entropy function, the maximum L*q*-likelihood estimator (ML*q*E). In our approach, the role of the observations is modified by slightly changing the model of reference by means of the distortion parameter *q*. From this standpoint, L*q*-likelihood estimation can be regarded as the minimization of the discrepancy between a distribution in a family and one that modifies the true distribution to diminish (or emphasize) the role of extreme observations.

In this framework, we provide theoretical insights concerning the statistical usage of the generalized entropy function. In particular, we highlight the role of the distortion parameter *q* and give the conditions that guarantee asymptotic efficiency of the ML*q*E. Further, the new methodology is shown to be very useful when estimating high-dimensional parameters and small tail probabilities. This aspect is important in many applications where we must deal with the fact that the number of observations available is not large in relation to the number of parameters or the probability of occurrence of the event of interest. Standard large sample theory guarantees that the maximum likelihood estimator (MLE) is asymptotically efficient, meaning that when the sample size is large, the MLE is at least as accurate as any other estimator. However, for a moderate or small sample size, it turns out that the ML$q$E can offer a dramatic improvement in terms of mean squared error at the expense of a slightly increased bias, as will be seen in our numerical results.

For finite sample performance of ML$q$E, not only the size of $q_n - 1$ but also its sign (i.e., the direction of distortion) is important. It turns out that for different families or different parametric functions of the same family, the beneficial direction of distortion can be different. In addition, for some parameters, ML$q$E does not produce any improvement. We have found that an asymptotic variance expression of the ML$q$E is very helpful to decide the direction of distortion for applications.

The paper is organized as follows. In Section 2, we examine some information-theoretical quantities and introduce the ML$q$E; in Section 3, we present its basic asymptotic properties for exponential families. In particular, a necessary and sufficient condition on the choice of $q$ in terms of the sample size to ensure a proper asymptotic normality and efficiency is established. A generalization that goes out of the exponential family is presented in Section 4. In Section 5, we consider the plug-in approach for tail probability estimation based on ML$q$E. The asymptotic properties of the plug-in estimator are derived and its efficiency is compared to the traditional MLE. In Section 6, we discuss the choice of the distortion parameter $q$. In Section 7, we present Monte Carlo simulations and examine the behavior of ML$q$E in finite sample situations. In Section 8, concluding remarks are given. Technical proofs of the theorems are deferred to Appendix A.

## 2. Generalized entropy and the maximum L$q$-likelihood estimator.

Consider a $\sigma$-finite measure $\mu$ on a measurable space $(\Omega, \mathscr{F})$. The Kullback–Leibler (KL) divergence [21, 22] (or relative entropy) between two density functions $g$ and $f$ with respect to $\mu$ is

$$(2.1) \qquad \mathcal{D}(f\|g) = E_f \log \frac{f(X)}{g(X)} = \int_\Omega f(x) \log \frac{f(x)}{g(x)} \, d\mu(x).$$

Note that finding the density $g$ that minimizes $\mathcal{D}(f\|g)$ is equivalent to minimizing Shannon's entropy [28] $\mathcal{H}(f, g) = -E_f \log g(X)$.

DEFINITION 2.1. Let $f$ and $g$ be two density functions. The $q$-*entropy* of $g$ with respect to $f$ is defined as

$$(2.2) \qquad \mathcal{H}_q(f, g) = -E_f L_q\{g(X)\}, \qquad q > 0,$$

where $L_q(u) = \log u$ if $q = 1$ and $L_q(u) = (u^{1-q} - 1)/(1-q)$, otherwise.

The function $L_q$ represents a Box–Cox transformation in statistics and in other contexts it is often called a deformed logarithm. Note that if $q \to 1$, then $L_q(u) \to \log(u)$ and the usual definition of Shannon's entropy is recovered.

Let $\mathscr{M} = \{f(x; \theta), \theta \in \Theta\}$ be a family of parametrized density functions and suppose that the true density of observations, denoted by $f(x; \theta_0)$, is a member of $\mathscr{M}$. Assume further that $\mathscr{M}$ is closed under the transformation

$$(2.3) \qquad f(x; \theta)^{(r)} = \frac{f(x; \theta)^r}{\int_\Omega f(x; \theta)^r \, d\mu(x)}, \qquad r > 0.$$

The transformed density $f(x; \theta)^{(r)}$ is often referred to as *zooming* or *escort* distribution [1, 7, 26] and the parameter $r$ provides a tool to accentuate different regions of the untransformed true density $f(x; \theta)$. In particular, when $r < 1$, regions with density values close to zero are accentuated, while for $r > 1$, regions with density values further from zero are emphasized.

Consider the following KL divergence between $f(x; \theta)$ and $f(x; \theta_0)^{(r)}$:

$$(2.4) \qquad \mathcal{D}_r(\theta_0 \| \theta) = \int_\Omega f(x; \theta_0)^{(r)} \log \frac{f(x; \theta_0)^{(r)}}{f(x; \theta)} \, d\mu(x).$$

Let $\theta^*$ be the value such that $f(x; \theta^*) = f(x; \theta_0)^{(r)}$ and assume that differentiation can be passed under the integral sign. Then, clearly $\theta^*$ minimizes $\mathcal{D}_r(\theta_0 \| \theta)$ over $\theta$. Let $\theta^{**}$ be the value such that $f(x; \theta^{**}) = f(x; \theta_0)^{(1/q)}$, $q > 0$. Since we have $\nabla_\theta \mathcal{H}_q(\theta_0, \theta)|_{\theta^{**}} = 0$ and $\nabla_\theta^2 \mathcal{H}_q(\theta_0, \theta)|_{\theta^{**}}$ is positive definite, $\mathcal{H}_q(\theta_0, \theta)$ has a minimum at $\theta^{**}$.

The derivations above show the minimizer of $\mathcal{D}_r(\theta_0 \| \theta)$ over $\theta$ is the same as the minimizer of $\mathcal{H}_r(\theta_0, \theta)$ over $\theta$ when $q = 1/r$. Clearly, by considering the divergence with respect to a distorted version of the true density we introduce a certain amount of bias. Nevertheless, the bias can be properly controlled by an adequate choice of the distortion parameter $q$, and later we shall discuss the benefits gained from paying such a price for parameter estimation. The next definition introduces the estimator based on the empirical version of the $q$-entropy.

DEFINITION 2.2. Let $X_1, \ldots, X_n$ be an i.i.d. sample from $f(x; \theta_0)$, $\theta_0 \in \Theta$. The maximum L$q$-likelihood estimator (ML$q$E) of $\theta_0$ is defined as

$$(2.5) \qquad \widetilde{\theta}_n = \arg\max_{\theta \in \Theta} \sum_{i=1}^n L_q[f(X_i; \theta)], \qquad q > 0.$$

When $q \to 1$, if the estimator $\widetilde{\theta}_n$ exists, then it approaches the maximum likelihood estimate of the parameters, which maximizes $\sum_i \log f(X_i; \theta)$. In this sense, the ML$q$E extends the classic method, resulting in a general inferential procedure that inherits most of the desirable features of traditional maximum likelihood, and at the same time can improve over MLE due to variance reduction, as will be seen.

Define

$$(2.6) \qquad \begin{aligned} U(x; \theta) &= \nabla_\theta \log\{f(x; \theta)\}, \\ U^*(X; \theta, q) &= U(X; \theta) f(X; \theta)^{1-q}. \end{aligned}$$

In general, the estimating equations have the form

$$(2.7) \qquad \sum_{i=1}^n U^*(X_i; \theta, q) = 0.$$

Equation (2.7) offers a natural interpretation of the ML$q$E as a solution to a *weighted* likelihood. When $q \neq 1$, (2.7) provides a relative-to-the-model re-weighting. Observations that disagree with the model receive low or high weight depending on $q < 1$ or $q > 1$. In the case $q = 1$, all the observations receive the same weight.

The strategy of setting weights that are proportional to a power transformation of the assumed density has some connections with the methods proposed by Windham [33], Basu et al. [6] and Choi, Hall and Presnell [8]. In these approaches, however, the main objective is robust estimation and the weights are set based on a fixed constant not depending on the sample size.

EXAMPLE 2.1. The simple but illuminating case of an exponential distribution will be used as a recurrent example in the course of the paper. Consider an i.i.d. sample of size $n$ from a distribution with density $\lambda_0 \exp\{-x\lambda_0\}$, $x > 0$ and $\lambda_0 > 0$. In this case, the $L_q$-likelihood equation is

$$(2.8) \qquad \sum_{i=1}^n e^{-[X_i \lambda - \log \lambda](1-q)} \left( -X_i + \frac{1}{\lambda} \right) = 0.$$

With $q = 1$, the usual maximum likelihood estimator is $\widehat{\lambda} = (\sum_i X_i / n)^{-1} = \overline{X}^{-1}$. However, when $q \neq 1$, (2.8) can be rewritten as

$$(2.9) \qquad \lambda = \left( \frac{\sum_{i=1}^n X_i w_i(X_i, \lambda, q)}{\sum_{i=1}^n w_i(X_i, \lambda, q)} \right)^{-1},$$

where $w_i := e^{-[X_i \lambda - \log \lambda](1-q)}$. When $q < 1$, the role played by observations corresponding to higher density values are accentuated; when $q > 1$, observations corresponding to density values close to zero are accentuated.

**3. Asymptotics of the ML*q*E for exponential families.** In this section, we present the asymptotic properties of the new estimator when the degree of distortion is chosen according to the sample size. In the remainder of the paper, we focus on exponential families, although some generalization results are presented in Section 4. In particular, we consider density functions of the form

$$f(x; \theta) = \exp\{\theta^{\mathsf{T}} b(x) - A(\theta)\}, \tag{3.1}$$

where $\theta \in \Theta \subseteq \mathbb{R}^p$ is a real valued natural parameter vector, $b(x)$ is the vector of functions with elements $b_j(x)$ $(j = 1, \ldots, p)$ and $A(\theta) = \log \int_{\Omega} e^{\theta^{\mathsf{T}} b(x)} d\mu(x)$ is the cumulant generating function (or log normalizer). For simplicity in presentation, the family is assumed to be of full rank (but similar results hold for curved exponential families). The true parameter will be denoted by $\theta_0$.

3.1. *Consistency.* Consider $\theta_n^*$, the value such that

$$E_{\theta_0} U^*(X; \theta_n^*, q_n) = 0. \tag{3.2}$$

It can be easily shown that $\theta_n^* = \theta_0/q_n$. Since the actual target of $\widetilde{\theta}_n$ is $\theta_n^*$, to retrieve asymptotic unbiasedness of $\widetilde{\theta}_n$, $q_n$ must converge to 1. We call $\theta_n^*$ the *surrogate* parameter of $\theta_0$. We impose the following conditions:

A.1 $q_n > 0$ is a sequence such that $q_n \to 1$ as $n \to \infty$.
A.2 The parameter space $\Theta$ is compact and the parameter $\theta_0$ is an interior point in $\Theta$.

In similar contexts, the compactness condition on $\Theta$ is used for technical reasons (see, e.g., Wang, van Eeden and Zidek [32]), as is the case here.

THEOREM 3.1. *Under assumptions* A.1 *and* A.2*, with probability going to* 1*, the* $L_q$*-likelihood equation yields a unique solution* $\widehat{\theta}_n$ *that is the maximizer of the* $L_q$*-likelihood function in* $\Theta$*. Furthermore, we have* $\widetilde{\theta}_n \xrightarrow{P} \theta_0$*.*

REMARK. When $\Theta$ is compact, the ML*q*E always exists under our conditions, although it is not necessarily unique with probability one.

3.2. *Asymptotic normality.*

THEOREM 3.2. *If assumptions* A.1 *and* A.2 *hold, then we have*

$$\sqrt{n} V_n^{-1/2}(\widetilde{\theta}_n - \theta_n^*) \xrightarrow{\mathscr{D}} N_p(0, \mathbf{I}_p) \qquad as\ n \to \infty, \tag{3.3}$$

*where* $\mathbf{I}_p$ *is the* $(p \times p)$ *identity matrix,* $V_n = J_n^{-1} K_n J_n^{-1}$ *and*

$$K_n = E_{\theta_0}[U^*(X; \theta_n^*, q_n)]^{\mathsf{T}}[U^*(X; \theta_n^*, q_n)], \tag{3.4}$$

$$J_n = E_{\theta_0}[\nabla_{\theta} U^*(X; \theta_n^*, q_n)]. \tag{3.5}$$

*A necessary and sufficient condition for asymptotic normality of MLqE around $\theta_0$ is $\sqrt{n}(q_n - 1) \to 0$.*

Let $m(\theta) := \nabla_\theta A(\theta)$ and $D(\theta) := \nabla_\theta^2 A(\theta)$. Note that $K_n$ and $J_n$ can be expressed as

$$(3.6) \qquad K_n = c_{2,n}(D(\theta_{2,n}) + [m(\theta_{2,n}) - m(\theta_n^*)][m(\theta_{2,n}) - m(\theta_n^*)]^{\mathsf{T}})$$

and

$$(3.7) \qquad \begin{aligned} J_n = {}& c_{1,n}(1 - q_n)D(\theta_{1,n}) - c_{1,n}D(\theta_n^*) \\ & + c_{1,n}(1 - q_n)[m(\theta_{1,n}) - m(\theta_n^*)][m(\theta_{1,n}) - m(\theta_n^*)]^{\mathsf{T}}, \end{aligned}$$

where $c_{k,n} = \exp\{A(\theta_{n,k}) - A(\theta_0)\}$ and $\theta_{k,n} = k\theta_0(1/q_n - 1) + \theta_0$. When $q_n \to 1$, it is seen that $V_n \to -D(\theta_0)$, the asymptotic variance of the MLE. When $\Theta \subseteq \mathbb{R}^1$ we use the notation $\sigma_n^2$ for the asymptotic variance in place of $V_n$. Note that the existence of moments are ensured by the functional form of the exponential families (e.g., see [23]).

REMARKS. (i) When $q$ is fixed, the MLqE is a regular M-estimator [18], which converges in probability to $\theta^* = \theta_0/q$. (ii) With the explicit expression of $\theta_n^*$, one may consider correcting the bias of MLqE by using the estimator $q_n\widetilde{\theta}_n$. The numerical results are not promising in this direction under correct model specification.

EXAMPLE 3.1 (Exponential distribution). The surrogate parameter is $\theta_n^* = \lambda_0/q_n$ and a lengthy but straightforward calculation shows that the asymptotic variance of the MLqE of $\lambda_0$ is

$$(3.8) \qquad \sigma_n^2 = \left(\frac{\lambda_0}{q_n}\right)^2 \left[\frac{q_n^2 - 2q_n + 2}{q_n^3(2 - q_n)^3}\right] \to \lambda_0^2$$

as $n \to \infty$. By Theorem 3.2, we conclude that $n^{1/2}\sigma_n^{-1}(\widetilde{\lambda}_n - \lambda_0/q_n)$ converges weakly to a standard normal distribution as $n \to \infty$. Clearly, the asymptotic calculation does not produce any advantage of MLqE in terms of reducing the limiting variance. However, for an interval of $q_n$, we have $\sigma_n^2 < \lambda_0^2$ (see Section 6) and, based on our simulations, an improvement of the accuracy is achieved in finite sample sizes as long as $0 < q_n - 1 = o(n^{-1/2})$, which ensures a proper asymptotic normality of $\widetilde{\lambda}_n$. For the re-scaled estimator $q_n\widehat{\lambda}_n$, the expression $q_n^2\sigma_n^2$ is larger than 1 unless $q = 1$, which suggests that $q_n\widehat{\lambda}_n$ may be at best no better than $\widehat{\lambda}_n$.

EXAMPLE 3.2 (Multivariate normal distribution). Consider a multivariate normal family with mean vector $\boldsymbol{\mu}$ and covariance matrix $\boldsymbol{\Sigma}$. Two convenient matrix operators in this setting are the vec($\cdot$) (vector) and vech($\cdot$) (vector-half). Namely, $\text{vec}\colon\mathbb{R}^{r\times p}\mapsto\mathbb{R}^{rp}$ stacks the columns of the argument matrix. For symmetric matrices, $\text{vech}\colon\mathbb{S}^{p\times p}\mapsto\mathbb{R}^{p(p+1)/2}$ stacks only the unique part of each column that lies on or below the diagonal [25]. Further, for a symmetric matrix $M$, define the extension matrix $G$ as $\text{vec}\,M = G\,\text{vech}\,M$. Thus, $\theta_0 = (\boldsymbol{\mu}^{\mathsf{T}}, \text{vech}^{\mathsf{T}}\,\boldsymbol{\Sigma})^{\mathsf{T}}$ and under such a parametrization, it is easy to show the surrogate parameter solving (3.2) is $\theta_n^* = (\boldsymbol{\mu}^{\mathsf{T}}, \sqrt{q_n}\,\text{vech}^{\mathsf{T}}\,\boldsymbol{\Sigma})^{\mathsf{T}}$, where interestingly the mean component does not depend on $q_n$. In fact, for symmetric distributions about the mean, it can be shown that the distortion imposed to the model affects the spread of the distribution but leaves the mean unchanged. Consequently, the ML$q$E is expected to influence the estimation of $\boldsymbol{\Sigma}$ without much effect on $\boldsymbol{\mu}$. This will be clearly seen in our simulation results (see Section 7.3). The calculation in Appendix B shows that the asymptotic variance of the ML$q$E of $\theta_0$ is the block-diagonal matrix

$$(3.9) \qquad V_n = \begin{pmatrix} \dfrac{(2-q)^{2+p}}{(3-2q)^{1+p/2}}\boldsymbol{\Sigma} & 0 \\[2ex] 0 & \begin{array}{c} \dfrac{4q^2[(3-2q)^2+1](2-q)^{4+p}}{[(2-q)^2+1]^2(3-2q)^{2+p/2}} \\ \times\,[G^{\mathsf{T}}(\boldsymbol{\Sigma}^{-1}\otimes\boldsymbol{\Sigma}^{-1})G]^{-1} \end{array} \end{pmatrix},$$

where $\otimes$ denotes the Kronecker product.

## 4. A generalization.

In this section, we relax the restriction of exponential family and present consistency and asymptotic normality results for ML$q$E under some regularity conditions.

THEOREM 4.1. *Let $q_n$ be a sequence such that $q_n \to 1$ as $n \to \infty$ and assume the following:*

B.1 $\theta_0$ *is an interior point in $\Theta$.*

B.2 $E_{\theta_0}\sup_{\theta\in\Theta}\|U(X;\theta)\|^2 < \infty$ *and* $E_{\theta_0}\sup_{\theta\in\Theta}[f(X;\theta)^\delta - 1]^2 \to 0$ *as* $\delta \to 0$.

B.3 $\sup_{\theta\in\Theta}\|\frac{1}{n}\sum_{i=1}^n U(X_i;\theta) - E_{\theta_0}U(X;\theta)\| \overset{p}{\to} 0$ *as* $n \to \infty$,

*where $\|\cdot\|$ denotes the $\ell_2$-norm. Then, with probability going to 1, the $L_q$-likelihood equation yields a unique solution $\widetilde{\theta}_n$ that maximizes the $L_q$-likelihood. Furthermore, we have $\widetilde{\theta}_n \overset{P}{\to} \theta_0$.*

REMARK 4.2. (i) Although for a large $n$ the $L_q$-likelihood equation has a unique zero with a high probability, for finite samples there may be roots that are actually bad estimates. (ii) The uniform convergence in condition B.3 is satisfied if the set of functions $\{U(x,\theta)\colon\theta\in\Theta\}$ is Glivenko–Cantelli under the true parameter $\theta_0$ (see, e.g., [31], Chapter 19.2). In particular, it suffices to require (i) $U(x;\theta)$ is continuous in $\theta$ for every $x$ and dominated by an integrable function and (ii) compactness of $\Theta$.

For each $\theta \in \Theta$, define a symmetric $p \times p$ matrix $I^*(x; \theta, q) = \nabla_\theta U^*(x; \theta, q)$, where $U^*$ represents the modified score function as in (2.6) and let the matrices $K_n$, $J_n$ and $V_n$ be as defined in the previous section.

THEOREM 4.3. *Let $q_n$ be a sequence such that $q_n \to 1$ and $\theta_n^* \to \theta_0$ as $n \to \infty$, where $\theta_n^*$ is the solution of $EU^*(X; \theta_n^*, q_n) = 0$. Suppose $U^*(x; \theta, q)$ is twice differentiable in $\theta$ for every $x$ and assume the following:*

C.1 $\max_{1 \le k \le p} E_{\theta_0} |U_k^*(X, \theta_n^*, q_n)|^3$, $k = 1, \dots, p$, *is upper bounded by a constant.*

C.2 *The smallest eigenvalue of $K_n$ is bounded away from zero.*

C.3 $E_{\theta_0} \{I^*(X, \theta_n^*, q_n)\}_{kl}^2$, $k, l = 1, \dots, p$, *are upper bounded by a constant.*

C.4 *The second-order partial derivatives of $U^*(x, \theta, q_n)$ are dominated by an integrable function with respect to the true distribution of $X$ for all $\theta$ in a neighborhood of $\theta_0$ and $q_n$ in a neighborhood of $1$.*

*Then,*

$$(4.1) \qquad \sqrt{n} V_n^{-1/2} (\widetilde{\theta}_n - \theta_n^*) \xrightarrow{\mathscr{D}} N_p(0, \mathbf{I}) \qquad as \ n \to \infty.$$

## 5. Estimation of the tail probability.
In this section, we address the problem of tail probability estimation, using the popular plug-in procedure, where the point estimate of the unknown parameter is substituted into the parametric function of interest. We focus on a one-dimensional case, that is, $p = 1$, and derive the asymptotic distribution of the plug-in estimator for the tail probability based on the MLq method. For an application of the MLqE proposed in this work on financial risk estimation, see Ferrari and Paterlini [12].

Let $\alpha(x; \theta) = P_\theta(X \le x)$ or $\alpha(x; \theta) = 1 - P_\theta(X \le x)$, depending on whether we are considering the lower tail or the upper tail of the distribution. Without loss of generality, we focus on the latter from now on, and assume $\alpha(x; \theta) > 0$ for all $x$ [of course $\alpha(x; \theta) \to 0$ as $x \to \infty$]. When $x$ is fixed, under some conditions, the familiar delta method shows that an asymptotically normally distributed and efficient estimator of $\theta$ makes the plug-in estimator of $\alpha(x; \theta)$ also asymptotically normal and efficient. However, in most applications a large sample size is demanded in order for this asymptotic behavior to be accurate for a small tail probability. As a consequence, the setup with $x$ fixed but $n \to \infty$ presents an overly optimistic view, as it ignores the possible difficulty due to smallness of the tail probability in relation to the sample size $n$. Instead, allowing $x$ to increase in $n$ (so that the tail probability to be estimated becomes smaller as the sample size increases) more realistically addresses the problem.

5.1. *Asymptotic normality of the plug-in MLq estimator.* We are interested in estimating $\alpha(x_n; \theta_0)$, where $x_n \to \infty$ as $n \to \infty$. For $\theta^* \in \Theta$ and $\delta > 0$, define

$$(5.1) \qquad \beta(x; \theta^*; \delta) = \sup_{\theta \in \Theta \cap [\theta^* - \delta/\sqrt{n}, \theta^* + \delta/\sqrt{n}]} \left| \frac{\alpha''(x; \theta)}{\alpha''(x; \theta^*)} \right|$$

and $\gamma(x; \theta) = \alpha''(x; \theta)/\alpha'(x; \theta)$.

THEOREM 5.1. *Let $\theta_n^*$ be as in the previous section. Under assumptions A.1 and A.2, if $n^{-1/2}|\gamma(x_n; \theta_n^*)|\beta(x_n; \theta_n^*; \delta) \to 0$ for each $\delta > 0$, then*

$$\sqrt{n} \frac{\alpha(x_n; \widetilde{\theta}_n) - \alpha(x_n; \theta_n^*)}{\sigma_n \alpha'(x_n; \theta_n^*)} \xrightarrow{\mathscr{D}} N(0, 1),$$

*where $\sigma_n = -[E_{\theta_0} U^*(X; \theta_n^*)^2]^{1/2}/E_{\theta_0}[\partial U^*(X; \theta, q_n)/\partial \theta|_{\theta_n^*}]$.*

REMARKS. (i) For the main requirement of the theorem on the order of the sequence $x_n$, it is easiest to be verified on a case by case basis. For instance, in the case of the exponential distribution in (A.4), for $x_n > 0$,

$$\beta(x_n; \lambda_n^*; \delta) = \sup_{\lambda \in \lambda_n^* \pm \delta/\sqrt{n}} \frac{e^{-x_n \lambda} x_n^2}{e^{-x_n \lambda_n^*} x_n^2} \leq \sup_{\lambda \in \lambda_n^* \pm \delta/\sqrt{n}} e^{x_n |\lambda - \lambda_n^*|} = e^{\delta x_n/\sqrt{n}}.$$

Moreover, $\gamma(x_n; \lambda_n^*) = -x_n$. So, the condition reads $n^{-1/2} x_n e^{\delta x_n/\sqrt{n}} \to 0$, that is, $n^{-1/2} x_n \to 0$. (ii) The plug-in estimator based on $q_n \widetilde{\theta}_n$ has been examined as well. With $q_n \to 1$, we did not find any significant advantage.

The condition $n^{-1/2}|\gamma(x_n; \theta_n^*)|\beta(x_n; \theta_n^*; \delta) \to 0$, to some degree, describes the interplay between the sample size $n$, $x_n$ and $q_n$ for the asymptotic normality to hold. When $x_n \to \infty$ too fast so as to violate the condition, the asymptotic normality is not guaranteed, which indicates the extreme difficulty in estimating a tiny tail probability. In the next section, we will use this framework to compare the MLqE of the tail probability, $\alpha(x_n; \widetilde{\theta}_n)$, with the one based on the traditional MLE, $\alpha(x_n; \widehat{\theta}_n)$.

In many applications, the quantity of interest is quantile instead of the tail probability. In our setting, the quantile function is defined as $\rho(s; \theta) = \alpha^{-1}(s; \theta)$, $0 < s < 1$ and $\theta \in \Theta$. Next, we present the analogue of Theorem 5.1 for the plug-in estimator of the quantile. Define

$$(5.2) \qquad \beta_1(s; \theta^*; \delta) = \sup_{\theta \in \Theta \cap [\theta^* - \delta/\sqrt{n}, \theta^* + \delta/\sqrt{n}]} \left| \frac{\rho''(s; \theta)}{\rho'(s; \theta^*)} \right|, \qquad \delta > 0,$$

and $\gamma_1(s; \theta) = \rho''(s; \theta)/\rho'(s; \theta)$.

THEOREM 5.2. *Let $0 < s_n < 1$ be a nonincreasing sequence such that $s_n \searrow 0$ as $n \to \infty$ and let $\theta_n^*$ and $q_n$ be as in Theorem 5.1. Under assumptions A.1 and A.2, for a sequence $s_n$ such that $n^{-1/2}|\gamma_1(s_n; \theta_n^*)||\beta_1(s_n; \theta_n^*; \delta) \to 0$ for each $\delta > 0$, we have*

$$\sqrt{n} \frac{\rho(s_n; \widetilde{\theta}_n) - \rho(s_n; \theta_n^*)}{\sigma_n \rho'(s_n; \theta_n^*)} \xrightarrow{\mathscr{D}} N(0, 1).$$

5.2. *Relative efficiency between MLE and MLqE.* In Section 3, we showed that when $(q_n - 1)\sqrt{n} \to 0$, the MLqE is asymptotically as efficient as the MLE. For tail probability estimation, with $x_n \to \infty$, it is unclear if the MLqE performs efficiently.

Consider $w_n$ and $v_n$, two estimators of a parametric function $g_n(\theta)$ such that both $\sqrt{n}(w_n - a_n)/\sigma_n$ and $\sqrt{n}(v_n - b_n)/\tau_n$ converge weakly to a standard normal distribution as $n \to \infty$ for some deterministic sequences $a_n$, $b_n$, $\sigma_n > 0$ and $\tau_n > 0$.

DEFINITION 5.1. Define

(5.3) $$\Lambda(w_n, v_n) := \frac{(b_n - g_n(\theta))^2 + \tau_n^2/n}{(a_n - g_n(\theta))^2 + \sigma_n^2/n}.$$

The *bias adjusted asymptotic relative efficiency* of $w_n$ with respect to $v_n$ is $\lim_{n \to \infty} \Lambda(w_n, v_n)$, provided that the limit exists.

It can be easily verified that the definition does not depend on the specific choice of $a_n$, $b_n$, $\sigma_n$ and $\tau_n$ among equivalent expressions.

COROLLARY 5.3. *Under the conditions of Theorem 5.1, when $q_n$ is chosen such that*

(5.4) $$n^{1/2}\alpha(x_n; \theta_n^*)\alpha(x_n; \theta_0)^{-1} \to 1 \quad and \quad \alpha'(x_n; \theta_n^*)\alpha'(x_n; \theta_0)^{-1} \to 1,$$

*then $\Lambda(\alpha(x_n; \widehat{\theta}_n), \alpha(x_n; \widetilde{\theta}_n)) = 1$.*

The result, which follows directly from Theorem 5.1 and Definition 5.1, says that when $q_n$ is chosen sufficiently close to 1, asymptotically speaking, the MLqE is as efficient as the MLE.

EXAMPLE 5.1 (Continued). In this case, we have $\alpha(x_n; \lambda) = e^{-\lambda x_n}$ and $\alpha'(x_n; \lambda) = -x_n e^{-\lambda x_n}$. For sequences $x_n$ and $q_n$ such that $x_n/\sqrt{n} \to 0$ and $(q_n - 1)\sqrt{n} \to 0$, we have that

(5.5) $$\sqrt{n} \frac{(e^{-\widetilde{\lambda}_n x_n} - e^{-\lambda_0/q_n x_n})}{\lambda_0 x_n e^{-\lambda_0/q_n x_n}} \xrightarrow{\mathscr{D}} N(0, 1).$$

When $q_n = 1$ for all $n$, we recover the usual plug-in estimator based on MLE. With the asymptotic expressions given above,

$$(5.6) \quad \Lambda(\alpha(x; \widehat{\lambda}_n), \alpha(x; \widetilde{\lambda}_n)) = \frac{n}{\lambda_0^2 x_n^2} (e^{-x_n(\lambda_0/q_n - \lambda_0)} - 1)^2 + e^{-2x_n(\lambda_0/q_n - \lambda_0)},$$

which is greater than 1 when $q_n > 1$. Thus, no advantage in terms of MSE is expected by considering $q_n > 1$ (which introduces bias and enlarges the variance at the same time).

Although in limits ML$q$E is not more efficient than MLE, ML$q$E can be much better than MLE due to variance reduction as will be clearly seen in Section 7. The following calculation provides a heuristic understanding. Let $r_n = 1 - 1/q_n$. Add and subtract 1 in (5.6), obtaining

$$(5.7) \quad \frac{nr_n^2 L_{1/q_n}(e^{-x_n \lambda_0})^2}{\lambda_0^2 x_n^2} + r_n L_{1/q_n}(e^{2x_n \lambda_0}) + 1 < nr_n^2 + r_n 2 x_n \lambda_0 + 1,$$

where the last inequality holds as $L_{1/q_n}(u) < \log(u)$ for any $u > 0$ and $q < 1$. Next, we impose (5.7) to be smaller than 1 and solve for $q_n$, obtaining

$$(5.8) \quad T_n := \left(1 + \frac{2\lambda_0 x_n}{n}\right)^{-1} < q_n < 1.$$

This provides some insights on the choice of the sequence $q_n$ in accordance to the size of the probability to be estimated. If $q_n$ approaches 1 too quickly from below, the gain obtained in terms of variance vanishes rapidly as $n$ becomes larger. On the other hand, if $q_n$ converges to 1 too slowly, the bias dominates the variance and the MLE outperforms the ML$q$E. This understanding is confirmed in our simulation study.

**6. On the choice of $q$.** For the exponential distribution example, we have observed the following:

1. For estimating the natural parameter, when $q_n \to 1$, the asymptotic variance of ML$q$E is equivalent to that of MLE in limit, but can be smaller. For instance, in the variance expression (3.8) one can easily check that $(q^2 - 2q + 2)/[q^5(2 - q)^3] < 1$ for $1 < q < 1.40$; thus, choosing the distortion parameter in such a range gives $\sigma_n^2 < \lambda_0^2$.

2. For estimating the tail probability, when $q_n \to 1$, the asymptotic variance of ML$q$E can be of a smaller order than that of MLE, although there is a bias that approaches 0. In particular:

    (i) ML$q$E cannot be asymptotically more efficient than MLE.

    (ii) ML$q$E is asymptotically as efficient as MLE when $q_n$ is chosen to be close enough to 1. In the case of tail probability for the exponential distribution, it suffices to choose $q_n$ such that $(q_n - 1)x_n \to 0$.

3. One approach to choosing $q$ is to minimize an estimated asymptotic mean squared error of the estimator when it is mathematically tractable. In the case of the exponential distribution, by Theorem 5.1 we have the following expression for the asymptotic mean squared error:

$$\text{MSE}(q, \lambda_0) = (e^{-\lambda_0/qx_n} - e^{-\lambda_0 x_n})^2$$

(6.1)

$$+ \left(\frac{\lambda_0}{q}\right)^2 n^{-1} \left(\frac{q^2 - 2q + 2}{q^3(2-q)^3}\right) x_n^2 e^{-2\lambda_0/qx_n}.$$

However, since $\lambda_0$ is unknown, we consider

(6.2)
$$q^* = \underset{q \in (0,1)}{\arg\min} \{\text{MSE}(q, \widehat{\lambda})\},$$

where $\widehat{\lambda}$ is the MLE. This will be also used in some of our simulation studies.

In general, unlike in the above example, closed-form expressions of the asymptotic mean squared error are not available, which calls for more work on this issue. In the literature on applications of nonextensive entropy, although some discussions on choosing $q$ have been made often from physical considerations, it is unclear how to do it from a statistical perspective. In particular, the direction of distortion (i.e., $q > 1$ or $q < 1$) needs to be decided. We offer the following observations and thoughts:

1. For estimating the parameters in an exponential family, although $|q_n - 1| n^{-1/2}$ guarantees the right asymptotic normality (i.e., asymptotic normality centered around $\theta_0$), one direction of distortion typically reduces the variance of estimation and consequently improves the MSE. In the exponential distribution case, $q_n$ needs to be slightly greater than 1, but for estimating the covariance matrix for multivariate normal observations, based on the asymptotic variance formula in Example 3.2, $q_n$ needs to be slightly smaller than 1. For a given family, the expression of the asymptotic covariance matrix for the MLqE given in Section 3 can be used to find the beneficial direction of distortion. Our numerical investigations confirm this understanding.

2. To minimize the mean squared error for tail probability estimation for the exponential distribution family, we need $0 < q_n < 1$. This choice is in the opposite direction for estimating the parameter $\lambda$ itself. Thus, the optimal choice of $q_n$ is not a characteristic of the family alone but also depends on the parametric function to be estimated.

3. For some parametric functions, the MLqE makes little change. For the multivariate normal family, the surrogate value of the mean parameter stays exactly the same while the variance parameters are altered.

4. We have found empirically that given the right distortion direction, choices of $q_n$ with $|1 - q_n|$ between $1/n$ and $1/\sqrt{n}$ usually improves—to different extents—over the MLE.

**7. Monte Carlo results.** In this section, the performance of the ML$q$E in finite samples is explored via simulations. Our study includes (i) an assessment of the accuracy for tail probability estimation and reliability of confidence intervals and (ii) an assessment of the performance of ML$q$E for estimating multidimensional parameters, including regression settings with generalized linear models. The standard MLE is used as a benchmark throughout the study.

In this section, we present both deterministic and data-driven approaches on choosing $q_n$. First, deterministic choices are used to explore the possible advantage of the ML$q$E for tail probability estimation with $q_n$ approaching 1 fast when $x$ is fixed and $q_n$ approaching 1 slowly when $x$ increases with $n$. Then, the data-driven choice in Section 6 is applied. For multivariate normal and GLM families, where estimation of the MSE or prediction error becomes analytically cumbersome, we choose $q_n = 1 - 1/n$, which satisfies $1 - q_n = o(n^{-1/2})$ that is needed for asymptotic normality around $\theta_0$. In all considered cases, numerical solution of (2.7) is found using variable metric algorithm (e.g., see Broyden [15]), where the ML solution is chosen as the starting value.

7.1. *Mean squared error: role of the distortion parameter $q$.* In the first group of simulations, we compare the estimators of the true tail probability $\alpha = \alpha(x; \lambda_0)$, obtained via the ML$q$ method and the traditional maximum likelihood approach. Particularly, we are interested in assessing the relative performance of the two estimators for different choices of the sample size by taking the ratio between the two mean squared errors, $\mathrm{MSE}(\hat{\alpha}_n)/\mathrm{MSE}(\tilde{\alpha}_n)$. The simulations are structured as follows: (i) For any given sample size $n \geq 2$, a number $B = 10{,}000$ of Monte Carlo samples $X_1, \ldots, X_n$ is generated from an exponential distribution with parameter $\lambda_0 = 1$. (ii) For each sample, the ML$q$ and ML estimates of $\alpha$, respectively, $\tilde{\alpha}_{n,k} = \alpha(x; \tilde{\lambda}_{n,k})$ and $\hat{\alpha}_{n,k} = \alpha(x; \hat{\lambda}_{n,k})$, $k = 1, \ldots, B$, are obtained. (iii) For each sample size $n$, the relative performance between the two estimators is evaluated by the ratio $\hat{R}_n = \mathrm{MSE_{MC}}(\hat{\alpha}_n)/\mathrm{MSE_{MC}}(\tilde{\alpha}_n)$, where $\mathrm{MSE_{MC}}$ denotes the Monte Carlo estimate of the mean squared error. In addition, let $\overline{y}_1 = B^{-1} \sum_{k=1}^{B} (\hat{\alpha}_{n,k} - \alpha)^2$ and $\overline{y}_2 = B^{-1} \sum_{k=1}^{B} (\tilde{\alpha}_{n,k} - \alpha)^2$. By the central limit theorem, for large values of $B$, $\overline{y} = (\overline{y}_1, \overline{y}_2)'$ approximately has a bi-variate normal distribution with mean $(\mathrm{MSE}(\hat{\alpha}_n), \mathrm{MSE}(\tilde{\alpha}_n))'$ and a certain covariance matrix $\Gamma$. Thus, the standard error for $\hat{R}_n$ can be computed by the delta method [11] as

$$\mathrm{se}(\hat{R}_n) = B^{-1/2} \left( \frac{\hat{\gamma}_{11}}{\overline{y}_2^2} - 2\hat{\gamma}_{12} \frac{\overline{y}_1}{\overline{y}_2^3} + \hat{\gamma}_{22} \frac{\overline{y}_1^2}{\overline{y}_2^4} \right)^{1/2},$$

where $\hat{\gamma}_{11}$, $\hat{\gamma}_{22}$ and $\hat{\gamma}_{12}$ denote, respectively, the Monte Carlo estimates for the components of the covariance matrix $\Gamma$.

*Case* 1: *fixed* $\alpha$ *and* $q$. Figure 1 illustrates the behavior of $\widehat{R}_n$ for several choices of the sample size. In general, we observe that for relatively small sample sizes, $\widehat{R}_n > 1$ and the MLqE clearly outperforms the traditional MLE. Such a behavior is much more accentuated for smaller values of the tail probability to be estimated. In contrast, when the sample size is larger, the bias component plays an increasingly relevant role and eventually we observe that $\widehat{R}_n < 1$. This case is presented in Figure 1(a) for values of the true tail probability $\alpha = 0.01, 0.005, 0.003$ and a fixed distortion parameter $q = 0.5$. Moreover, the results presented in Figure 1(b) show that smaller values of the distortion parameter $q$ accentuate the benefits attainable in a small sample situation.

*Case* 2: *fixed* $\alpha$ *and* $q_n \nearrow 1$. In the second experimental setting, illustrated in Figure 2(a), the tail probability $\alpha$ is fixed, while we let $q_n$ be a sequence such that $q_n \nearrow 1$ and $0 < q_n < 1$. For illustrative purposes we choose the sequence $q_n = [1/2 + e^{0.3(n-20)}]/[1 + e^{0.3(n-20)}]$, $n \geq 2$, and study $R_n$ for different choices of the true tail probability to be estimated. For small values of the sample size, the chosen sequence $q_n$ converges relatively slowly to 1 and the distortion parameter produces benefits in terms of variance. In contrast, when the sample size becomes larger, $q_n$ adjusts quickly to one. As a consequence, for large samples the MLqE exhibits the same behavior shown by the traditional MLE.

*Case* 3: $\alpha_n \searrow 0$ *and* $q_n \nearrow 1$. The last experimental setting of this subsection examines the case where both the true tail probability and the distortion parameter change depending on the sample size. We consider sequences of distortion parameters converging slowly relative to the sequence of quantiles $x_n$. In particular we set $q_n = 1 - [10 \log(n + 10)]^{-1}$ and $x_n = n^{1/(2+\delta)}$. In the simulation described in Figure 2(b), we illustrate the behavior of the estimator for $\delta = 0.5, 1.0$ and $1.5$, confirming the theoretical findings discussed in Section 5.

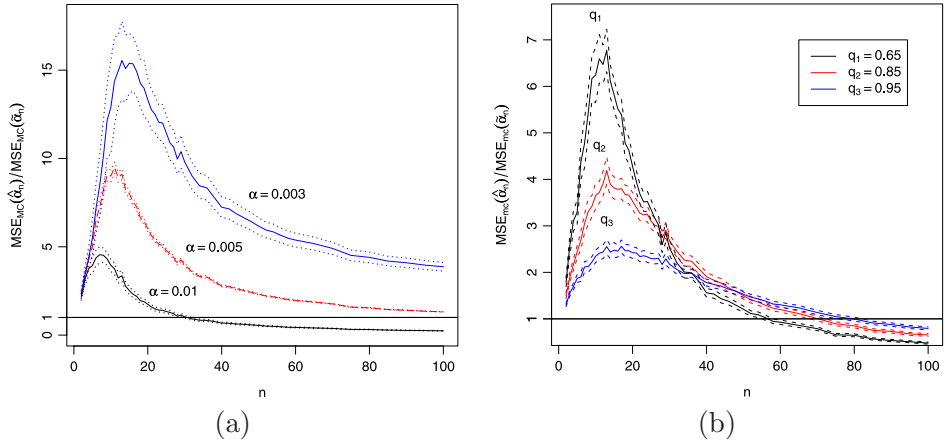

(a)                                                (b)

FIG. 1. *Monte Carlo mean squared error ratio computed from* $B = 10{,}000$ *samples of size* $n$. *In* (a) *we use a fixed distortion parameter* $q = 0.5$ *and true tail probability* $\alpha = 0.01, 0.005, 0.003$. *The dashed lines represent* 99% *confidence bands. In* (b) *we set* $\alpha = 0.003$ *and use* $q_1 = 0.65$, $q_2 = 0.85$ *and* $q_3 = 0.95$. *The dashed lines represent* 90% *confidence bands.*

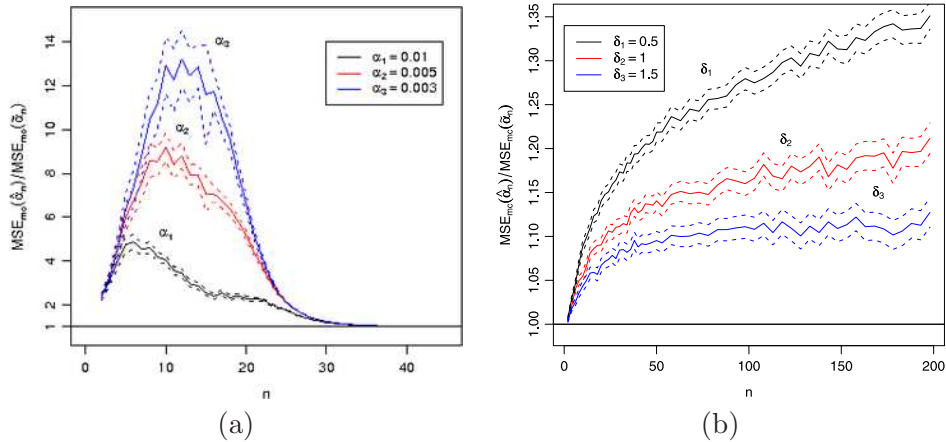

Fig. 2. (a) *Monte Carlo mean squared error ratio computed from $B = 10,000$ samples of size $n$, for different values of the true probability: $\alpha_1 = 0.01$, $\alpha_2 = 0.005$ and $\alpha_3 = 0.003$. The distortion parameter is computed as $q_n = [1/2 + e^{0.3(n-20)}]/[1 + e^{0.3(n-20)}]$. (b) Monte Carlo mean squared error ratio computed from $B = 10,000$ samples of size $n$. We use sequences $q_n = 1 - [10 \log(n + 10)]^{-1}$ and $x_n = n^{1/(2+\delta)}$ ($\delta_1 = 0.5$, $\delta_2 = 1.0$ and $\delta_3 = 1.5$). The dashed lines represent 99% confidence bands.*

### 7.2. Asymptotic and bootstrap confidence intervals.

The main objective of the simulations presented in this subsection is twofold: (a) to study the reliability of MLqE based confidence intervals constructed using three commonly used methods: asymptotic normality, parametric bootstraps and nonparametric bootstraps; (b) to compare the results with those obtained using MLE. The structure of simulations is similar to that of Section 7.1, but a data-driven choice of $q_n$ is used. (i) For each sample, first we compute $\widehat{\lambda}_n$, the MLE of $\lambda_0$. We substitute $\widehat{\lambda}_n$ in (6.1) and solve it numerically in order to obtain $q^*$ as described there. (ii) For each sample, the MLq and ML estimates of the tail probability $\alpha$ are obtained. The standard errors of the estimates are computed using three different methods: the asymptotic formula derived in (5.5), nonparametric bootstrap and parametric bootstrap. The number of replicates employed in bootstrap re-sampling is 500. We construct 95% bootstrap confidence intervals based on the bootstrap quantiles and check the coverage of the true value $\alpha$.

*MC means and standard deviations of estimators of $\alpha$, along with the MC mean of the standard error computed using: (i) asymptotic normality, (ii) bootstrap and (iii) parametric bootstrap. The true tail probability is $\alpha = 0.01$ and $q = 1$ corresponds to the MLE*

| $n$ | $q^*$ | Estimate | St. dev. | $\text{se}_{\text{asy}}$ | $\text{se}_{\text{boot}}$ | $\text{se}_{\text{pboot}}$ |
|---|---|---|---|---|---|---|
| 15 | 0.939 | 0.009489 | 0.010975 | 0.010472 | 0.011923 | 0.010241 |
|  | 1.000 | 0.013464 | 0.014830 | 0.013313 | 0.013672 | 0.015090 |
| 25 | 0.959 | 0.009693 | 0.008417 | 0.008470 | 0.009134 | 0.008298 |
|  | 1.000 | 0.012108 | 0.010517 | 0.009919 | 0.010227 | 0.010950 |
| 50 | 0.977 | 0.010108 | 0.006261 | 0.006326 | 0.006575 | 0.006249 |
|  | 1.000 | 0.011385 | 0.007354 | 0.006894 | 0.007083 | 0.007318 |
| 100 | 0.988 | 0.010158 | 0.004480 | 0.004568 | 0.004680 | 0.004549 |
|  | 1.000 | 0.010789 | 0.004908 | 0.004778 | 0.004880 | 0.004943 |
| 500 | 0.998 | 0.010006 | 0.002014 | 0.002052 | 0.002061 | 0.002050 |
|  | 1.000 | 0.010122 | 0.002055 | 0.002070 | 0.002073 | 0.002087 |

In Table 1, we show the Monte Carlo means of $\widehat{\alpha}_n$ and $\widetilde{\alpha}_n$, their standard deviations and the standard errors computed with the three methods described above. In addition, we report the Monte Carlo average of the estimates of optimal distortion parameter $q^*$. When $q^* = 1$, the results refer to the MLE case. Not surprisingly, $q^*$ approaches 1 as the sample size increases. When the sample size is small, the ML$q$E has a smaller standard deviation and better performance. When $n$ is larger, the advantage of ML$q$E diminishes. As far as the standard errors are concerned, the asymptotic method and the parametric bootstrap seem to provide values somewhat closer to the Monte Carlo standard deviation for the considered sample sizes.

In Table 2, we compare the accuracy of 95% confidence intervals and report the relative length of the intervals for ML$q$E over those for MLE. Although the coverage probability for ML$q$E is slightly smaller than that of MLE (in the order of 1%), we observe a substantial reduction in the interval length for all of the considered cases. The most evident benefits occur when the sample size is small. Furthermore, in general, the intervals computed via parametric bootstrap outperform the other two methods in terms of coverage and length.

7.3. *Multivariate normal distribution.* In this subsection, we evaluate the ML$q$ methodology for estimating the mean and covariance matrix of a multivariate normal distribution. We generate $B = 10{,}000$ samples from a multivariate normal $N_p(\mu, \Sigma)$, where $\mu$ is the $p$-dimensional unknown mean vector and $\Sigma$ is the unknown $(p \times p)$ covariance matrix. In our simulation, the true mean is $\mu = 0$ and the $ij$th element of $\Sigma$ is $\rho^{|i-j|}$, where $-1 < \rho < 1$. To gauge performance for the mean we employed the usual $L_2$-norm. For

TABLE 2

*MC coverage rate of 95% confidence intervals for $\alpha$, computed using* (i) *asymptotic normality,* (ii) *boostrap and* (iii) *parametric bootstrap. RL is the length of the intervals of MLqE over that of MLE. The true tail probability is $\alpha = 0.01$ and $q = 1$ corresponds to the MLE*

| | | **Asympt.** | | **Boot.** | | **Par. boot.** | |
|---|---|---|---|---|---|---|---|
| $n$ | $q^*$ | Coverage (%) | RL | Coverage (%) | RL | Coverage (%) | RL |
| 15 | 0.939 | 79.2 | 0.787 | 89.1 | 0.865 | 92.9 | 0.657 |
| | 1.000 | 80.9 | | 88.4 | | 92.5 | |
| 25 | 0.958 | 83.4 | 0.854 | 91.8 | 0.890 | 93.6 | 0.733 |
| | 1.000 | 84.3 | | 90.8 | | 94.2 | |
| 50 | 0.977 | 87.1 | 0.918 | 92.3 | 0.928 | 93.9 | 0.824 |
| | 1.000 | 88.4 | | 91.6 | | 93.4 | |
| 100 | 0.988 | 91.1 | 0.956 | 93.3 | 0.960 | 94.7 | 0.889 |
| | 1.000 | 92.2 | | 92.9 | | 94.3 | |
| 500 | 0.998 | 94.5 | 0.991 | 95.0 | 0.995 | 95.2 | 0.962 |
| | 1.000 | 94.7 | | 94.6 | | 94.8 | |

the covariance matrix, we considered the loss function

$$\Delta(\Sigma, \widehat{\Sigma}_q) = \operatorname{tr}(\Sigma^{-1}\widehat{\Sigma}_q - \mathbf{I})^2, \tag{7.1}$$

where $\widehat{\Sigma}_q$ represents the ML$q$ estimate of $\Sigma$ with $q = 1 - 1/n$. Note that the loss is 0 when $\Sigma = \widehat{\Sigma}_q$ and is positive otherwise. Moreover, the loss is invariant to the transformations $A\Sigma A^{\mathsf{T}}$ and $A\widehat{\Sigma}_q A^{\mathsf{T}}$ for a nonsingular matrix $A$. The use of such a loss function is common in literature (e.g., Huang et al. [17]).

In Table 3, we show simulation results for moderate or small sample sizes ranging from 10 to 100 for various dimensions of the covariance matrix $\Sigma$. The entries in the table represent the Monte Carlo mean of $\Delta(\Sigma, \widehat{\Sigma}_1)$ over that of $\Delta(\Sigma, \widehat{\Sigma}_q)$, where $\widehat{\Sigma}_1$ is the usual ML estimate multiplied by the correction factor $n/(n-1)$. The standard error of the ratio is computed via the delta method. Clearly, the ML$q$E performs well for smaller sample sizes. Interestingly, the squared error for the ML$q$E reduces dramatically compared to that of the MLE as the dimension increases. Remarkably, when $p = 8$ the gain in accuracy persists even for larger sample sizes, ranging from about 22% to 84%. We tried various structures of $\Sigma$ and obtained performances comparable to the ones presented. For $\mu$ we found that ML$q$E performs nearly identically to MLE for all choices of $p$ and $n$, which is not surprising given the findings in Section 3. For brevity we omit the results on $\mu$.

7.4. *Generalized linear models.* Our methodology can be promptly extended to the popular framework of the generalized linear models. Con-sider the regression setting where each outcome of the dependent variables, $Y$, is drawn from a distribution in the exponential family. The mean $\eta$ of the distribution is assumed to depend on the independent variables, $X$, through $E(Y|X) = \eta = g^{-1}(X^{\mathsf{T}}\beta)$, where $X$ is the design matrix, $\beta$ is a $p$-dimensional vector of unknown parameters and $g$ is the link function. In our simulations, we consider two notable instances: (i) $Y$ from an exponential distribution with $\eta = \exp(-x^{\mathsf{T}}\beta)$; (ii) $Y$ from a Bernoulli distribution with $\eta = 1/(1 + \exp\{x^{\mathsf{T}}\beta\})$. The first case represents the exponential regression model, which is a basic setup for time-to-event analysis. The latter is the popular logistic regression model.

TABLE 3
*Monte Carlo mean of $\Delta(\Sigma, \widehat{\Sigma}_1)$ over that of $\Delta(\Sigma, \widehat{\Sigma}_q)$ with standard error in parenthesis*

| | $p$ | | | |
|---|---|---|---|---|
| $n$ | **1** | **2** | **4** | **8** |
| 10 | 1.225 (0.018) | 1.298 (0.019) | 1.740 (0.029) | 1.804 (0.022) |
| 15 | 1.147 (0.014) | 1.249 (0.017) | 1.506 (0.021) | 1.840 (0.026) |
| 25 | 1.083 (0.011) | 1.153 (0.012) | 1.313 (0.016) | 1.562 (0.020) |
| 50 | 1.041 (0.007) | 1.052 (0.007) | 1.199 (0.011) | 1.377 (0.015) |
| 100 | 1.018 (0.005) | 1.033 (0.005) | 1.051 (0.006) | 1.222 (0.011) |

We initialize the simulations by generating design points randomly drawn from the unit hypercube $[-1, 1]^p$. The entries of the true vector of coefficients $\beta$ are assigned by sampling $p$ points at random in the interval $[-1, 1]$, obtaining values $\beta = (-0.57, 0.94, 0.16, -0.72, 0.68, 0.92, 0.80, 0.04, 0.64, 0.34, 0.38, 0.47)$. The values of $X$ and $\beta$ are kept fixed during the simulations. Then, 1000 Monte Carlo samples of $Y|X$ are generated according to the two models described above and for each sample ML$q$ and ML estimates are computed. The prediction error based on independent out-of-sample observations is

$$(7.2) \qquad \mathrm{PE}_q = \frac{1}{10^3} \sum_{j=1}^{10^3} (Y_j^{\mathrm{test}} - g^{-1}(X_j^{\mathrm{test}} \widehat{\beta}_{q,n}))^2,$$

where $\widehat{\beta}_{q,n}$ is the ML$q$E of $\beta$. In Table 4 we present the prediction error for various choices of $n$ and $p$. For both models, the ML$q$E outperforms the classic MLE for all considered cases. The benefits from ML$q$E can be remarkable when the dimension of the parameter space is larger. This is particularly evident in the case of the exponential regression, where the prediction error of MLE is at least twice that of ML$q$E. In one case, when $n = 25$ and $p = 12$, the ML$q$E is about nine times more accurate. This is mainly due to ML$q$E's stabilization of the variance component, which for the MLE tends to become large quickly when $n$ is very small compared to $p$.

TABLE 4
*Monte Carlo mean of* $\mathrm{PE}_1$ *over that of* $\mathrm{PE}_q$ *for exponential and logistic regression with standard error in parenthesis*

| | $n$ | | | |
|---|---|---|---|---|
| $p$ | **25** | **50** | **100** | **250** |
| | Exp. regression | | | |
| 2 | 2.549 (0.003) | 2.410 (0.002) | 2.500 (0.003) | 2.534 (0.003) |
| 4 | 2.469 (0.002) | 2.392 (0.002) | 2.543 (0.002) | 2.493 (0.002) |
| 8 | 4.262 (0.012) | 2.941 (0.004) | 3.547 (0.006) | 3.582 (0.006) |
| 12 | 9.295 (0.120) | 3.644 (0.008) | 3.322 (0.005) | 5.259 (0.027) |
| | Logistic regression | | | |
| 2 | 1.156 (0.006) | 1.329 (0.006) | 1.205 (0.003) | 1.385 (0.003) |
| 4 | 1.484 (0.022) | 1.141 (0.003) | 1.502 (0.007) | 1.353 (0.003) |
| 8 | 1.178 (0.008) | 1.132 (0.003) | 1.290 (0.004) | 1.300 (0.002) |
| 12 | 1.086 (0.005) | 1.141 (0.003) | 1.227 (0.003) | 1.329 (0.002) |

Although for the logistic regression we observe a similar behavior, the gain in high dimension becomes more evident for larger $n$.

**8. Concluding remarks.** In this work, we have introduced the ML$q$E, a new parametric estimator inspired by a class of generalized information measures that have been successfully used in several scientific disciplines. The ML$q$E may also be viewed as a natural extension of the classical MLE. It can preserve the large sample properties of the MLE, while—by means of a distortion parameter $q$—allowing modification of the trade-off between bias and variance in small or moderate sample situations. The Monte Carlo simulations support that when the sample size is small or moderate, the ML$q$E can successfully trade bias for variance, obtaining a reduction of the mean squared error, sometimes very dramatically.

Overall, this work makes a significant contribution to parametric estimation and applications of nonextensive entropies. For parametric models, MLE is by far the most commonly used estimator and the substantial improvement as seen in our numerical work seems relevant and important to applications. Given the increasing attention to $q$-entropy in other closely related disciplines, our theoretical results provide a useful view from a statistical perspective. For instance, from the literature, although $q$ is chosen from interesting physical considerations, for statistical estimation (e.g., financial data analysis where $q$-entropy is considered), there are few clues as to how to choose the direction and amount of distortion.

Besides the theoretical optimality results and often remarkably improved performance over MLE, our proposed method is very practical in terms of implementability and computational efficiency. The estimating equations are simply obtained by replacing the logarithm of log-likelihood function in the usual maximum likelihood procedure by the distorted logarithm. Thus, the resulting optimization task can be easily formulated in terms of a weighted version of the familiar score function, with weights proportional to the $(1 - q)$th power of the assumed density. Hence, similarly to other techniques based on re-weighing of the likelihood, simple and fast algorithms for solving the ML$q$ equations numerically (possibly even for large problems) can be derived.

For the ML$q$ estimators, helpful insights on their behaviors may be gained from robust analysis. For a given $q$, (2.6) defines an M-estimator of the surrogate parameter $\theta^*$. It seems that global robustness properties, such as a high breakdown point, may be established for a properly chosen distortion parameter, which would add value to the ML$q$ methodology.

High-dimensional estimation has recently become a central theme in statistics. The results in this work suggest that the ML$q$ methodology may be a valuable tool for some high-dimensional estimation problems (such as gamma regression and covariance matrix estimation as demonstrated in this paper) as a powerful remedy to the MLE. We believe this is an interesting direction for further exploration.

Finally, more research on the practical choices of $q$ and their theoretical properties will be valuable. To this end, higher-order asymptotic treatment of the distribution (or moments) of the ML$q$E will be helpful. For instance, derivation of saddle-point approximations of order $n^{-3/2}$, along the lines of Field and Ronchetti [13] and Daniels [10], may be profitably used to give improved approximations of the MSE.

## APPENDIX A: PROOFS

In all of the following proofs we denote $\psi_n(\theta) := n^{-1} \sum_{i=1}^{n} \nabla_\theta L_{q_n}(f(X_i; \theta))$. For exponential families, since $f(x; \theta) = e^{\theta^\mathsf{T} b(x) - A(\theta)}$, we have

$$(A.1) \qquad \psi_n(\theta) = \frac{1}{n} \sum_{i=1}^{n} e^{(1-q_n)(\theta^\mathsf{T} b(X_i) - A(\theta))} (b(X_i) - m(\theta)),$$

where $m(\theta) = \nabla_\theta A(\theta)$. The ML$q$ equation sets $\psi_n(\theta) = 0$ and solves for $\theta$. Moreover, we define $\varphi(x, \theta) := \theta^\mathsf{T} b(x) - A(\theta)$, and thus $f(x; \theta) = e^{\varphi(x,\theta)}$. When clear from the context, $\varphi(x, \theta)$ is denoted by $\varphi$.

**Proof of Theorem 3.1.** Define $\psi(\theta) := E_{\theta_0} \nabla_\theta \log (f(X; \theta))$. Since $f$ has the form in (3.1), we can write $\psi(\theta) = E_{\theta_0}[b(X) - m(\theta)]$. We want to show uniform convergence of $\psi_n(\theta)$ to $\psi(\theta)$ for all $\theta \in \Theta$ in probability. Clearly,

$$\sup_{\theta \in \Theta} \left\| \frac{1}{n} \sum_{i=1}^{n} e^{(1-q_n)(\theta^\mathsf{T} b(X_i) - A(\theta))} (b(X_i) - m(\theta)) - E_{\theta_0}[b(X) - m(\theta)] \right\|_1$$

$$
\begin{aligned}
\text{(A.2)} \quad &\leq \sup_{\theta \in \Theta} \left\| \frac{1}{n} \sum_{i=1}^{n} (e^{(1-q_n)(\theta^{\mathsf{T}} b(X_i) - A(\theta))} - 1)(b(X_i) - m(\theta)) \right\|_1 \\
&\quad + \sup_{\theta \in \Theta} \left\| \frac{1}{n} \sum_{i=1}^{n} (b(X_i) - m(\theta)) - E_{\theta_0}[b(X) - m(\theta)] \right\|_1,
\end{aligned}
$$

where $\|\cdot\|_1$ denotes the $\ell_1$-norm. Note that the second summand in (A.2) actually does not depend on $\theta$ [as $m(\theta)$ cancels out] and it converges to zero in probability by the law of large numbers. Next, let $s(X_i; \theta) := e^{(1-q_n)(\theta^{\mathsf{T}} b(X_i) - A(\theta))} - 1$ and $t(X_i; \theta) := b(X_i) - m(\theta)$. By the Cauchy–Schwarz inequality, the first summand in (A.2) is upper bounded by

$$
\text{(A.3)} \quad \sup_{\theta \in \Theta} \left\{ \sum_{j=1}^{p} \sqrt{\frac{1}{n} \sum_{i=1}^{n} s(X_i; \theta)^2} \sqrt{\frac{1}{n} \sum_{i=1}^{n} t_j(X_i; \theta)^2} \right\},
$$

where $t_j$ denotes the $j$th element of the vector $t(X_i; \theta)$. It follows that for (A.2), it suffices to show $n^{-1} \sum_i \sup_\theta s(X_i; \theta)^2 \xrightarrow{p} 0$ and $n^{-1} \sum_i \sup_\theta t_j(X_i; \theta)^2$ is bounded in probability. Since $\Theta$ is compact, $\sup_\theta |m(\theta)| \leq (c_1, c_1, \ldots, c_1)$ for some positive constant $c_1 < \infty$, and we have

$$
\text{(A.4)} \quad \frac{1}{n} \sum_{i=1}^{n} \sup_{\theta \in \Theta} t_j(X_i; \theta)^2 \leq \frac{2}{n} \sum_{i=1}^{n} b_j(X_i)^2 + 2(c_1)^2,
$$

where the last inequality from the basic fact that $(a - b)^2 \leq 2a^2 + 2b^2$ $(a, b \in \mathbb{R})$. The last expression in (A.4) is bounded in probability by some constant and $E_{\theta_0} b_j(X)^2 < \infty$ for all $j = 1, \ldots, p$. Next, note that

$$
\begin{aligned}
\text{(A.5)} \quad \sup_{\theta \in \Theta} \frac{1}{n} \sum_{i=1}^{n} s(X_i; \theta)^2 &\leq \frac{1}{n} \sum_{i=1}^{n} \sup_{\theta \in \Theta} e^{2(1-q_n)(\theta^{\mathsf{T}} b(X_i) - A(\theta))} \\
&\quad - \frac{2}{n} \sum_{i=1}^{n} \inf_{\theta \in \Theta} e^{(1-q_n)(\theta^{\mathsf{T}} b(X_i) - A(\theta))} + 1.
\end{aligned}
$$

Thus, to show $n^{-1} \sum_i \sup_\theta s(X_i; \theta)^2 \xrightarrow{p} 0$, it suffices to obtain $n^{-1} \times \sum_i \sup_\theta e^{2(1-q_n)\varphi(\theta)} - 1 \xrightarrow{p} 0$ and $n^{-1} \sum_i \inf_\theta e^{(1-q_n)\varphi(\theta)} - 1 \xrightarrow{p} 0$. Actually, since $\Theta$ is compact and $\sup_\theta e^{-A(\theta)} < c_2$ for some $c_2 < \infty$,

$$
\text{(A.6)} \quad \frac{1}{n} \sum_{i=1}^{n} \sup_{\theta \in \Theta} e^{2(1-q_n)(\theta^{\mathsf{T}} b(X_i) - A(\theta))} \leq \frac{1}{n} \sum_{i=1}^{n} e^{2|1-q_n|(|\log c_2| + \theta^{\mathsf{T}} |b(X_i)|)},
$$

where $\theta_j^{(*)} = \max\{|\theta_{j,0}^{(*)}|, |\theta_{j,1}^{(*)}|\}$, $j = 1, \ldots, p$ and $(\theta_{j,0}^{(*)}, \theta_{j,1}^{(*)})$ represent elementwise boundary points of $\theta_j$. For $r = 1, 2$,

$$
\text{(A.7)} \quad E_{\theta_0}[e^{2|1-q_n|(|\log c_2| + \theta^{(*)\mathsf{T}} |b(X)|)}]^r
$$

$$(A.8) \qquad = e^{2r|1-q_n||\log c_2| - A(\theta_0)} \int e^{[2r|1-q_n|\operatorname{sign}\{b(x)\}\theta^{(*)} + \theta_0]^{\mathsf{T}} b(x)} \, d\mu(x).$$

We decompose $\Omega$ into $2^p$ subsets in terms of the sign of the elements of $b(x)$. That is, $\Omega = \bigcup_{k=1}^{2^p} B_k$, where

$$B_1 = \{x \in \Omega : b_1(x) \geq 0, b_2(x) \geq 0, \ldots, b_{p-1}(x) \geq 0, b_p(x) \geq 0\},$$

$$(A.9) \quad B_2 = \{x \in \Omega : b_1(x) \geq 0, b_2(x) \geq 0, \ldots, b_{p-1}(x) \geq 0, b_p(x) < 0\},$$

$$B_3 = \{x \in \Omega : b_1(x) \geq 0, b_2(x) \geq 0, \ldots, b_{p-1}(x) < 0, b_p(x) \geq 0\}$$

and so on. Note that $\operatorname{sign}\{b(x)\}$ stays the same for each $B_i$, $i = 1, \ldots, 2^p$. Also because $\theta_0$ is an interior point, when $|1 - q_n|$ is small enough, the integral in (A.7) on $B_i$ is finite and by dominated convergence theorem,

$$(A.10) \quad \int_{B_k} e^{[2r|1-q_n|\operatorname{sign}\{b(x)\}\theta^{(*)} + \theta_0]^{\mathsf{T}} b(x)} \, d\mu(x) \overset{n \to \infty}{\to} \int_{B_k} e^{\theta_0^{\mathsf{T}} b(x)} \, d\mu(x).$$

Consequently,

$$(A.11) \quad \int e^{[2r|1-q_n|\operatorname{sign}\{b(x)\}\theta^{(*)} + \theta_0]^{\mathsf{T}} b(x)} \, d\mu(x) \overset{n \to \infty}{\to} \int e^{\theta_0^{\mathsf{T}} b(x) - A(\theta_0)} \, d\mu(x) = 1.$$

It follows that the mean and the variance of $\sup_\theta e^{2(1-q_n)[\theta^{\mathsf{T}} b(X) - A(\theta)]}$ converge to 1 and 0, respectively, as $n \to \infty$. Therefore, a straightforward application of Chebyshev's inequality gives

$$(A.12) \qquad \frac{1}{n} \sum_{i=1}^{n} \sup_{\theta \in \Theta} e^{2(1-q_n)[\theta^{\mathsf{T}} b(X_i) - A(\theta)]} \overset{p}{\to} 1, \qquad n \to \infty.$$

An analogous argument shows that

$$(A.13) \qquad \frac{1}{n} \sum_{i=1}^{n} \inf_{\theta \in \Theta} e^{(1-q_n)[\theta^{\mathsf{T}} b(X_i) - A(\theta)]} \overset{p}{\to} 1, \qquad n \to \infty.$$

Therefore, we have established $n^{-1} \sum_i \sup_\theta s(X_i; \theta)^2 \overset{p}{\to} 0$. Hence, (A.2) converges to zero in probability. By applying Lemma 5.9 on page 46 in [31], we know that with probability converging to 1, the solution of the ML$q$ equations is unique and it maximizes the ML$q$E.

**Proof of Theorem 3.2.** By Taylor's theorem, there exist a random point $\widetilde{\theta}$, in the line segment between $\theta_n^*$ and $\widetilde{\theta}_n$, such that with probability converging to one we have

$$(A.14) \quad \begin{aligned} 0 &= \psi_n(\underline{X}; \widetilde{\theta}_n) \\ &= \psi_n(\underline{X}; \theta_n^*) + \dot{\psi}_n(\underline{X}; \theta_n^*)(\widetilde{\theta}_n - \theta_n^*) + \tfrac{1}{2}(\widetilde{\theta}_n - \theta_n^*)^{\mathsf{T}} \ddot{\psi}_n(\underline{X}; \widetilde{\theta})(\widetilde{\theta}_n - \theta_n^*), \end{aligned}$$

where $\dot{\psi}_n$ is a $p \times p$ matrix of first-order derivatives and, similarly to page 68 in van der Vaart [31], $\ddot{\psi}_n$ denotes a $p$-vector of $(p \times p)$ matrices of second-order derivatives, respectively; $\underline{X}$ denotes the data vector. We can rewrite the above expression as

$$
(A.15) \qquad -\sqrt{n}\dot{\psi}_n(\theta_n^*)^{-1}\psi_n(\underline{X};\theta_n^*)
$$

$$
= \dot{\psi}_n(\theta_n^*)^{-1}\dot{\psi}_n(\underline{X};\theta_n^*)\sqrt{n}(\widetilde{\theta}_n - \theta_n^*)
$$

$$
(A.16) \qquad + \dot{\psi}_n(\theta_n^*)^{-1}\frac{\sqrt{n}}{2}(\widetilde{\theta}_n - \theta_n^*)^{\mathsf{T}}\ddot{\psi}_n(\underline{X};\widetilde{\theta})(\widetilde{\theta}_n - \theta_n^*),
$$

where $\dot{\psi}(\theta) = E_{\theta_0}\nabla_\theta^2 L_{q_n}f(X;\theta)$. Note that

$$
(A.17) \qquad \dot{\psi}(\theta) = E_{\theta_0}e^{(1-q_n)\varphi(\theta)}[(1-q_n)\nabla_\theta\varphi(\theta)^{\mathsf{T}}\nabla_\theta\varphi(\theta) - \nabla_{\theta\theta}^2\varphi(\theta)]
$$

$$
(A.18) \qquad = K_{1,n}E_{\mu_{1,n}}[(1-q_n)\nabla_\theta\varphi(\theta)^{\mathsf{T}}\nabla_\theta\varphi(\theta) - \nabla_{\theta\theta}^2\varphi(\theta)],
$$

where $\mu_{k,n} = k(1-q_n)\theta + \theta_0$ and $K_{k,n} = e^{A(\mu_{n,k})-A(\theta_0)}$. For $k,l \in \{1,\ldots,p\}$, we have

$$
(A.19) \qquad \{E_{\mu_{n,1}}\nabla_\theta\varphi(\theta)^{\mathsf{T}}\nabla_\theta\varphi(\theta)\}_{kl}
$$

$$
= E_{\mu_{n,1}}[(b_k(X) - m_k(\theta))(b_l(X) - m_l(\theta))]
$$

$$
(A.20) \qquad = E_{\mu_{n,1}}[(b_k(X) - m_k(\mu_{n,1}) + m_k(\mu_{n,1}) - m_k(\theta))
$$

$$
(A.21) \qquad \times (b_l(X) - m_l(\mu_{n,1}) + m_l(\mu_{n,1}) - m_l(\theta))]
$$

$$
(A.22) \qquad = E_{\mu_{n,1}}[(b_k(X) - m_k(\mu_{n,1}))(b_l(X) - m_l(\mu_{n,1}))]
$$

$$
(A.23) \qquad + [(m_k(\mu_{n,1}) - m_k(\theta))(m_l(\mu_{n,1}) - m_l(\theta))],
$$

where the first term in the last passage is the $kl$th element of the covariance matrix $-D(\theta)$ evaluated at $\mu_{n,1}$. Since $\Theta$ is compact, $\{\dot{\psi}(\theta)\}_{kl} \leq C_{kl}^* < \infty$, for some constants $C_{kl}^*$, $k,l \in \{1,\ldots,p\}$. We take the following steps to derive asymptotic normality.

*Step* 1. We first show that the left-hand side of (A.16) converges in distribution. Define the vector $Z_{n,i} := \nabla_\theta L_{q_n}f(X_i,\theta_n^*) - E_{\theta_0}\nabla_\theta L_{q_n}f(X_i,\theta_n^*)$ in $\mathbb{R}^p$. Consider an arbitrary vector $a \in \mathbb{R}^p$ and let $W_{n,i} := a^{\mathsf{T}}Z_{n,i}$ and $\overline{W}_n = n^{-1}\sum_i W_{n,i}$. Since $W_{n,i}$ $(1 \leq i \leq n)$ form a triangular array where $W_{n,i}$ are rowwise i.i.d., we check the Lyapunov condition. In our case, the condition reads

$$
(A.24) \qquad n^{-1/3}(EW_{n,1}^2)^{-1}(E[W_{n,1}^3])^{2/3} \to 0 \qquad \text{as } n \to \infty.
$$

Next, denote $\mu_{n,k} = \theta_0 + k(1-q_n)\theta_n^*$. One can see that

$$
(E[W_{n,1}^3])^{1/3} = K_n\left(E_{\mu_{n,3}}\left[\sum_{j=1}^p a_j(b_j(X) - m_j(\theta_n^*))\right]^3\right)^{2/3},
$$

where $K_n = \exp\{-\frac{2}{3}A(\theta_0) - 2(1-q_n)A(\theta_n^*) + \frac{2}{3}A(\mu_{n,3})\}$ and $K_n \to 1$ as $n \to \infty$. Since $\theta_0$ is an interior point in $\Theta$ (compact) the above quantity is uniformly upper bounded in $n$ by some finite constant. Next, consider

$$E[W_{n,1}^2] = E[a^\mathsf{T} Z_{n,1} Z_{n,1}^\mathsf{T}] = a^\mathsf{T} E[Z_{n,1} Z_{n,1}^\mathsf{T}]a.$$

A calculation similar to that in (A.23) for the matrix $Z_{n,1} Z_{n,1}^\mathsf{T}$ shows that the above quantity satisfies

$$(A.25) \qquad a^\mathsf{T}[-D(\mu_{n,2}) + M_n]a \to -a^\mathsf{T} D(\theta_0)a > 0, \qquad n \to \infty,$$

where the $kl$th element of $M_n$ is

$$(A.26) \qquad \{M_n\}_{kl} = (m_k(\mu_{n,2}) - m_k(\theta_n^*))(m_l(\mu_{n,2}) - m_l(\theta_n^*))$$

and $\mu_{n,2} \to \theta_0$ and $\theta_n^* \to \theta_0$, as $n \to \infty$. This shows that condition (A.24) holds and $\sqrt{n}(E[W_{n,1}^2])^{-1/2}a^\mathsf{T}\overline{W}_n \overset{\mathscr{D}}{\to} N_1(0,1)$. Hence, by the Cramér–Wold device (e.g., see [31]), we have

$$(A.27) \qquad \sqrt{n}[EZ_{n,1} Z_{n,1}^\mathsf{T}]^{-1/2}\overline{W}_n \overset{\mathscr{D}}{\to} N_p(0, \mathbf{I}_p).$$

*Step* 2. Next, we want convergence in probability of $\dot{\psi}(\theta_n^*)^{-1}\dot{\psi}_n(\underline{X}, \theta_n^*)$ to $\mathbf{I}_p$. For $k, l \in \{1, \dots, p\}$, given $\varepsilon > 0$, we have

$$(A.28) \qquad \begin{aligned} &P_{\theta_0}(|\{\dot{\psi}_n(\underline{X}, \theta_n^*)\}_{kl} - \{\dot{\psi}(\theta_n^*)\}_{kl}| > \varepsilon) \\ &\qquad \le n^{-1}\varepsilon^{-2}E_{\theta_0}\left[\frac{\partial^2}{\partial\theta_k\theta_l}L_{q_n}(f(X;\theta))\Big|_{\theta_n^*}\right]^2 \end{aligned}$$

by the i.i.d. assumption and Chebyshev's inequality. When $|1 - q_n| \le 1$, the expectation in (A.28) is

$$\begin{aligned} &E_{\theta_0}[e^{2(1-q_n)\varphi(\theta_n^*)}[(1-q_n)(b_k(X) - m_k(\theta_n^*))(b_l(X) - m_l(\theta_n^*)) + D(\theta_n^*)^2]]^2 \\ &\le 2E_{\mu_{n,2}}[(b_k(X) - m_k(\theta_n^*))(b_l(X) - m(\theta_n^*))^2 + D(\theta_n^*)^4] \\ &\qquad \times \exp\{-A(\theta_0) - 2(1-q_n)A(\theta_n^*) + A(\mu_{n,2})\}, \end{aligned}$$

where the inequality passage follows from the triangle inequality. Since $\Theta$ is compact and the existence of fourth moments is ensured for exponential families, the above quantity is upper bounded by some finite constant. Therefore, the right-hand side of (A.28) is upper bounded by a constant that converges to zero as $n \to \infty$. Since convergence in probability holds for each $k, l \in \{1, \dots, p\}$ and $p < \infty$, we have that the matrix difference $|\dot{\psi}_n(\underline{X}, \theta_n^*) - \dot{\psi}(\theta_n^*)|$ converges in probability to the zero matrix. From the calculation carried out in (A.17), one can see that $\dot{\psi}(\theta_n^*)$ is a deterministic sequence such that $\dot{\psi}(\theta_n^*) \to \dot{\psi}(\theta_0) = -\nabla_\theta^2 A(\theta_0)$, as $n \to \infty$. Thus, we have

$$(A.29) \quad |\dot{\psi}_n(\underline{X}; \theta_n^*) - \dot{\psi}(\theta_0)| \le |\dot{\psi}_n(\underline{X}; \theta_n^*) - \dot{\psi}(\theta_n^*)| + |\dot{\psi}(\theta_n^*) - \dot{\psi}(\theta_0^*)| \overset{p}{\to} 0$$

as $n \to \infty$. Therefore, $\dot{\psi}(\theta_n^*)^{-1}\dot{\psi}_n(\underline{X}, \theta_n^*) \xrightarrow{p} \mathbf{I}_p$.

*Step* 3. Here, we show that the second term on the right-hand side of (A.16) is negligible. Let $g(\underline{X}; \theta)$ be an element of the array $\ddot{\psi}_n(\underline{X}, \theta)$ of dimension $p \times p \times p$. For some fixed $\overline{\overline{\theta}}$ in the line segment between $\widetilde{\theta}$ and $\theta_n^*$, we have that

$$(A.30) \quad |g(\underline{X}; \widetilde{\theta}) - g(\underline{X}; \theta_n^*)| = |\nabla_\theta g(\underline{X}, \overline{\overline{\theta}})^{\mathsf{T}}||\widetilde{\theta} - \theta_n^*| \leq \sup_{\theta \in \Theta} |\nabla_\theta g(\underline{X}, \theta)||\widetilde{\theta} - \theta_n^*|.$$

A calculation shows that the $h$th element of the gradient vector in the expression above is

$$
\begin{aligned}
(A.31) \quad & \{\nabla_\theta g(\underline{X}, \theta)\}_h \\
&= n^{-1} \sum_{i=1}^n e^{(1-q_n)\varphi(\theta)} [(1-q_n)^3 \varphi(\theta)^{(1)} + (1-q_n)^2 \varphi(\theta)^{(2)} \\
&\hspace{4cm} + (1-q_n)\varphi(\theta)^{(3)} + \varphi(\theta)^{(4)}]
\end{aligned}
$$

for $h \in \{1, \ldots, p\}$, where $\varphi^{(k)}$ denotes the product of the partial derivatives of order $k$ with respect to $\theta$. As shown before in the proof of Theorem 3.1, $\sup_\theta e^{(1-q_n)\varphi(X_i, \theta)}$ has finite expectation when $|1-q_n|$ is small enough. Thus, by Markov's inequality, $\sup_\theta |g'(\underline{X}, \theta)|$ is bounded in probability. In addition, recall that the deterministic sequence $\dot{\psi}(\theta_n^*)$ converges to a constant. Hence, $\dot{\psi}(\theta_n^*)^{-1}\ddot{\psi}_n(\underline{X}; \widetilde{\theta}_0)$ is bounded in probability.

Since the third term in the expansion (A.16) is of higher order than the second term, by combining steps 1, 2 and 3 and applying Slutsky's lemma we obtain the desired asymptotic normality result.

**Proof of Theorem 4.1.** Uniform convergence of $\psi_n(\theta)$ to $\psi(\theta)$ for all $\theta \in \Theta$ in probability is satisfied if

$$\sup_{\theta \in \Theta} \left\| \frac{1}{n} \sum_{i=1}^n f(X_i; \theta)^{1-q_n} U(X_i; \theta) - E_{\theta_0} U(X, \theta) \right\|_1 \xrightarrow{p} 0.$$

The left-hand side of the above expression is upper bounded by

$$
\begin{aligned}
(A.32) \quad & \sup_{\theta \in \Theta} \left\| \frac{1}{n} \sum_{i=1}^n (f(X_i; \theta)^{1-q_n} - 1) U(X_i; \theta) \right\|_1 \\
& + \sup_{\theta \in \Theta} \left\| \frac{1}{n} \sum_{i=1}^n U(X_i; \theta) - E_{\theta_0} U(X, \theta) \right\|_1.
\end{aligned}
$$

By the Cauchy–Schwarz inequality, the first term of the above expression is upper bounded by

$$\sup_{\theta \in \Theta} \left\{ \sum_{j=1}^{p} \sqrt{\frac{1}{n} \sum_{i=1}^{n} (f(X_i; \theta)^{1-q_n} - 1)^2} \sqrt{\frac{1}{n} \sum_{i=1}^{n} U_j(X_i; \theta)^2} \right\}.$$

By assumption B.2, $n^{-1} \sum_i \sup_\theta U_j(X_i; \theta)^2$ is bounded in probability. Moreover, given $\epsilon > 0$, by Markov's inequality we have

$$P\left( n^{-1} \sum_i \sup_\theta (f(X_i; \theta)^{1-q_n} - 1)^2 > \epsilon \right) \le \epsilon^{-1} E \sup_\theta (f(X; \theta)^{1-q_n} - 1)^2,$$

which converges to zero by assumption B.2. By assumption B.3, the second summand in (A.32) converges to zero in probability.

**Proof of Theorem 4.3.** By Taylor's theorem, for a solution of the ML$q$ equation, there exists a random point $\overline{\overline{\theta}}$ between $\widetilde{\theta}_n$ and $\theta_n^*$ such that

$$
\begin{aligned}
0 = \frac{1}{n} \sum_{i=1}^{n} U^*(X_i, \theta_n^*, q_n) &+ \frac{1}{n} \sum_{i=1}^{n} \nabla_\theta U^*(X_i, \theta_n^*, q_n)(\widetilde{\theta}_n - \theta_n^*) \\
&+ \frac{1}{2}(\widetilde{\theta}_n - \theta_n^*)^{\mathsf{T}} \frac{1}{n} \sum_{i=1}^{n} \nabla_\theta^2 U^*(X_i, \overline{\overline{\theta}}, q_n)(\widetilde{\theta}_n - \theta_n^*).
\end{aligned}
$$
(A.33)

From Theorem 4.1, we know that with probability approaching 1, $\widetilde{\theta}_n$ is the unique ML$q$E and the above equation holds. Define $Z_{n,i} := U^*(X_i; \theta_n^*, q_n)$, $i = 1, \ldots, n$, a triangular array of i.i.d. random vectors and let $a \in \mathbb{R}^p$ be a vector of constants. Let $W_{n,i} := a^{\mathsf{T}} Z_{n,i}$. The Lyapunov condition for ensuring asymptotic normality of the linear combination $a^{\mathsf{T}} \sum_{i=1}^{n} Z_{n,i}/n$ for $a \in \mathbb{R}^p$ and $\|a\| > 0$ in this case reads

$$n^{-1/3} (E W_{n,1}^2)^{-1} (E[W_{n,1}^3])^{2/3} \to 0 \qquad \text{as } n \to \infty.$$

Under C.1 and C.2, this can be easily checked. The Cramér–Wold device implies

$$C_n \frac{1}{n} \sum_{i=1}^{n} U^*(X_i, \theta_n^*, q_n) \xrightarrow{\mathscr{D}} N_p(0, \mathbf{I}_p),$$

where $C_n := \sqrt{n}[E_{\theta_0} U^*(X, \theta_n^*)^{\mathsf{T}} U^*(X, \theta_n^*)]^{-1/2}$.

Next, consider the second term in (A.33). Given $\epsilon > 0$, for $k, l \in \{1, \ldots, p\}$, by Chebyshev's inequality

$$P\left( \left| \left\{ n^{-1} \sum_{i=1}^{n} I^*(X_i, \theta_n^*, q_n) \right\}_{k,l} - \{J_n\}_{k,l} \right| > \epsilon \right) \le \epsilon^{-2} n^{-2} E\{I^*(X, \theta_n^*, q_n)\}_{k,l}^2.$$

Thus, the right-hand side of the above expression converges to zero as $n \to \infty$ under C.3. Since convergence in probability is ensured for each $k, l$ and $p < \infty$, under C.2, we have that $|n^{-1} \sum_i I^*(X_i, \theta_n^*) - J_n|$ converges to the zero matrix in probability.

Finally, $n^{-1} \nabla_\theta^2 \sum_{i=1}^n U^*(X_i, \overline{\overline{\theta}}, q_n)$ in the third term of the expansion (A.33) is a $p \times p \times p$ array of partial second-order derivatives. By assumption, there is a neighborhood $B$ of $\theta_0$ for which each entry of $\nabla_\theta^2 U^*(x, \theta, q_n)$ is dominated by $g_0(x)$ for some $g_0(x) \geq 0$ for all $\theta \in B$. With probability tending to 1,

$$\left\| n^{-1} \sum_{i=1}^n \nabla_\theta^2 U^*(X_i, \overline{\overline{\theta}}, q_n) \right\| \leq p^3 n^{-1} \sum_{i=1}^n |g_0(X_i)|,$$

which is bounded in probability by the law of large numbers. Since the third term in the expansion (A.33) is of higher order than the second term, the normality result follows by applying Slutsky's lemma.

**Proof of Theorem 5.1.** From the second-order Taylor expansion of $\alpha(x_n; \widetilde{\theta}_n)$ about $\theta_n^*$ one can obtain

$$\text{(A.34)} \quad \begin{aligned} &\sqrt{n} \frac{(\alpha(x_n; \widetilde{\theta}_n) - \alpha(x_n; \theta_n^*))}{\sigma_n \alpha'(x_n; \theta_n^*)} \\ &= \sqrt{n} \frac{(\widetilde{\theta}_n - \theta_n^*)}{\sigma_n} + \frac{1}{2\sigma_n} \frac{\alpha''(x_n; \widetilde{\theta})}{\alpha'(x_n; \theta_n^*)} \sqrt{n}(\widetilde{\theta}_n - \theta_n^*)^2 \\ &= \sqrt{n} \frac{(\widetilde{\theta}_n - \theta_n^*)}{\sigma_n} + \frac{1}{2\sigma_n} \frac{\alpha''(x_n; \theta_n^*)}{\alpha'(x_n; \theta_n^*)} \frac{\alpha''(x_n; \widetilde{\theta})}{\alpha''(x_n; \theta_n^*)} \sqrt{n}(\widetilde{\theta}_n - \theta_n^*)^2, \end{aligned}$$

where $\widetilde{\theta}$ is a value between $\widetilde{\theta}_n$ and $\theta_n^*$. We need to show that the second term in (A.34) converges to zero in probability, that is,

$$\text{(A.35)} \quad \frac{\alpha''(x_n; \theta_n^*)}{\alpha'(x_n; \theta_n^*)} \frac{\alpha''(x_n; \widetilde{\theta})}{\alpha''(x_n; \theta_n^*)} \frac{\sigma_n}{\sqrt{n}} \frac{n(\widetilde{\theta}_n - \theta_n^*)^2}{\sigma_n^2} \xrightarrow{p} 0.$$

Since $\sqrt{n}(\widetilde{\theta}_n - \theta_n^*)/\sigma_n \xrightarrow{\mathscr{D}} N(0, 1)$ and $\sigma_n$ is upper bounded, we need

$$\text{(A.36)} \quad \frac{\alpha''(x_n; \theta_n^*)}{\alpha'(x_n; \theta_n^*)\sqrt{n}} \frac{\alpha''(x_n; \widetilde{\theta})}{\alpha''(x_n; \theta_n^*)} \xrightarrow{p} 0.$$

This holds under the assumptions of the theorem. This completes the proof of the theorem.

**Proof of Theorem 5.2.** The rationale presented here is analogous to that of Theorem 5.1. From the second-order Taylor expansion of $\rho(\widetilde{\theta}_n, s)$ about $\theta_n^*$ one can obtain

$$
\begin{aligned}
\text{(A.37)} \quad & \sqrt{n} \frac{\rho(s_n; \widetilde{\theta}_n) - \rho(s_n; \theta_n^*)}{\sigma_n \rho'(s_n; \theta_n^*)} \\
&= \sqrt{n} \frac{(\widetilde{\theta}_n - \theta_n^*)}{\sigma_n} + \frac{1}{2\sigma_n} \frac{\rho''(s_n; \overline{\overline{\theta}})}{\rho'(s_n; \theta_n^*)} \sqrt{n}(\widetilde{\theta}_n - \theta_n^*)^2,
\end{aligned}
$$

where $\overline{\overline{\theta}}$ is a value between $\widetilde{\theta}_n$ and $\theta_n^*$. The assumptions combined with Theorem 3.2 imply that the second term in (A.37) converges to 0 in probability. Hence, the central limit theorem follows from Slutsky's lemma.

## APPENDIX B: MULTIVARIATE NORMAL $N_P(\mu, \Sigma)$. ASYMPTOTIC DISTRIBUTION OF THE MLQE OF $\Sigma$

The log-likelihood function of a multivariate normal is

$$
\text{(B.1)} \quad \ell(\theta) = \log f(\mathbf{x}; \boldsymbol{\mu}, \boldsymbol{\Sigma}) = -\frac{p}{2}(2\pi) - \frac{1}{2} \log |\boldsymbol{\Sigma}| - \frac{1}{2}(\mathbf{x} - \boldsymbol{\mu})^{\mathsf{T}} \boldsymbol{\Sigma} (\mathbf{x} - \boldsymbol{\mu}).
$$

Recall that the surrogate parameter is $\theta^* = (\boldsymbol{\mu}^{\mathsf{T}}, q \operatorname{vech}^{\mathsf{T}} \boldsymbol{\Sigma})^{\mathsf{T}}$. The asymptotic variance is computed as $V = J^{-1}(\theta^*) K(\theta^*) J^{-1}(\theta^*)$, where

$$
\text{(B.2)} \quad K(\theta^*) = E_{\theta_0}[f(\mathbf{x}; \theta^*)^{2(1-q)} U(\mathbf{x}; \theta^*)^{\mathsf{T}} U(\mathbf{x}; \theta^*)]
$$

$$
\text{(B.3)} \quad = c_2 E^{(2)}[U(\mathbf{x}; \theta^*)^{\mathsf{T}} U(\mathbf{x}; \theta^*)]
$$

and

$$
\text{(B.4)} \quad J(\theta^*) = -q E_{\theta_0}[f(\mathbf{x}; \theta^*)^{1-q} U(\mathbf{x}; \theta^*)^{\mathsf{T}} U(\mathbf{x}; \theta^*)]
$$

$$
\text{(B.5)} \quad = -q c_1 E^{(1)}[U(\mathbf{x}; \theta^*)^{\mathsf{T}} U(\mathbf{x}; \theta^*)],
$$

where $E^{(r)}$ denotes expectation taken with respect to a normal with mean $\boldsymbol{\mu}$ and covariance matrix $[r(1-q) + 1]^{-1} \boldsymbol{\Sigma}$, $r = 1, 2$, and the normalizing constant $c_r$ is

$$
\text{(B.6)} \quad c_r := E_{\theta_0}[f(\mathbf{x}; \theta^*)^{r(1-q)}] = \frac{\int e^{-(r(1-q)+1)/2 \cdot (\mathbf{x} - \boldsymbol{\mu})^{\mathsf{T}} \boldsymbol{\Sigma}^{-1} (\mathbf{x} - \boldsymbol{\mu})} \, d\mathbf{x}}{(2\pi)^{rp(1-q)/2} |q \boldsymbol{\Sigma}|^{r(1-q)/2} (2\pi)^{1/2} |\boldsymbol{\Sigma}|^{1/2}}
$$

$$
\text{(B.7)} \quad = \frac{(r(1-q)+1)^{-p/2}}{(2\pi q^p |\boldsymbol{\Sigma}|)^{r(1-q)/2}}.
$$

Note that $K$ and $J$ can be partitioned into block form

$$
\text{(B.8)} \quad K = \begin{bmatrix} K_{11} & K_{12} \\ K_{21} & K_{22} \end{bmatrix}, \qquad J = \begin{bmatrix} J_{11} & J_{12} \\ J_{21} & J_{22} \end{bmatrix},
$$

where $K_{11}$ and $J_{11}$ depend on second-order derivatives of $U$ with respect to $\boldsymbol{\mu}$, $K_{22}$ and $J_{22}$ depend on second-order derivatives with respect to vech $\boldsymbol{\Sigma}$. The off-diagonal matrices $K_{12}$, $K_{21}$ depend on mixed derivatives of $U$ with respect $\boldsymbol{\mu}$ and vech$^{\mathsf{T}}$ $\boldsymbol{\Sigma}$. Since the mixed moments of order three are zero, one can check that $K_{21} = K_{12}^{\mathsf{T}} = 0$. Consequently, only the calculation of $K_{11}, K_{22}, J_{11}$ and $J_{22}$ is required and the expression of the asymptotic variance is given by

$$(B.9) \qquad V = \begin{bmatrix} V_{11} & 0 \\ 0 & V_{22} \end{bmatrix} := \begin{bmatrix} J_{11}^{-1} K_{11} J_{11}^{-1} & 0 \\ 0 & J_{22}^{-1} K_{22} J_{22}^{-1} \end{bmatrix}.$$

Next, we compute the entries of $K$ and $J$ using the approach employed by McCulloch [25] for the usual log-likelihood function. First, we use standard matrix differentiation to compute $K_{11}$ and $J_{11}$,

$$(B.10) \qquad K_{11} = c_2 E^{(2)}[(q\boldsymbol{\Sigma})^{-1}(\mathbf{x} - \boldsymbol{\mu})^{\mathsf{T}}(\mathbf{x} - \boldsymbol{\mu})(q\boldsymbol{\Sigma})^{-1}]$$

$$(B.11) \qquad = c_2 q^{-2}[2(1-q)+1]^{-1}\boldsymbol{\Sigma}^{-1}$$

and similarly one can obtain $J_{11} = -c_1 q^{-1}[(1-q)+1]^{-1}\boldsymbol{\Sigma}^{-1}$. Some straightforward algebra gives

$$(B.12) \qquad V_{11} = J_{11}^{-1} K_{11} J_{11}^{-1} = \frac{(2-q)^{2+p}}{(3-2q)^{1+p/2}}\boldsymbol{\Sigma}.$$

Next, we compute $V_{22}$. Let $\mathbf{z} := \boldsymbol{\Sigma}^{-1/2}(\mathbf{x} - \boldsymbol{\mu})$ using the following relationship derived by McCulloch ([25], page 682):

$$E[\nabla_{\text{vech}\,\boldsymbol{\Sigma}}\ell(\theta)]^{\mathsf{T}}[\nabla_{\text{vech}\,\boldsymbol{\Sigma}}\ell(\theta)]$$

$$(B.13) \qquad = 1/4 G^{\mathsf{T}}(\boldsymbol{\Sigma}^{-1/2} \otimes \boldsymbol{\Sigma}^{-1/2})(E[(\mathbf{z} \otimes \mathbf{z})(\mathbf{z}^{\mathsf{T}} \otimes \mathbf{z}^{\mathsf{T}})] - \text{vec}\,\mathbf{I}_p\,\text{vec}^{\mathsf{T}}\,\mathbf{I}_p)$$

$$\times (\boldsymbol{\Sigma}^{-1/2} \otimes \boldsymbol{\Sigma}^{-1/2})G.$$

Moreover, a result by Magnus and Neudecker ([24], page 388) shows

$$(B.14) \qquad E[(\mathbf{z} \otimes \mathbf{z})(\mathbf{z}^{\mathsf{T}} \otimes \mathbf{z}^{\mathsf{T}})] = \mathbf{I}_p + K_{p,p} + \text{vec}\,\mathbf{I}_p\,\text{vec}^{\mathsf{T}}\,\mathbf{I}_p,$$

where $K_{p,p}$ denotes the commutation matrix (see Magnus and Neudecker [24]). To compute $K_{22}$ and $J_{22}$, we need to evaluate (B.13) at $\theta^* = (\boldsymbol{\mu}^{\mathsf{T}}, q\,\text{vech}^{\mathsf{T}}\,\boldsymbol{\Sigma})^{\mathsf{T}}$, replacing the expectation operator with $c_r E^{(r)}[\cdot]$. In particular,

$$\{E^{(r)}[(\mathbf{z} \otimes \mathbf{z})(\mathbf{z}^{\mathsf{T}} \otimes \mathbf{z}^{\mathsf{T}})]_{\theta^*} - \text{vec}\,\mathbf{I}_p\,\text{vec}^{\mathsf{T}}\,\mathbf{I}_p\}G$$

$$(B.15) \qquad = (r(1-q)+1)^{-2}\{\mathbf{I}_p + K_{p,p}\}G$$

$$= 2(r(1-q)+1)^{-2}G,$$

where the last equality follows from the fact that $K_{p,p}G = G$. Therefore,

$$(B.16) \quad K_{22} = 1/(4q^2)c_2 G^{\mathrm{T}}(\boldsymbol{\Sigma}^{-1/2} \otimes \boldsymbol{\Sigma}^{-1/2})$$

$$(B.17) \quad \times (E[(\mathbf{z} \otimes \mathbf{z})(\mathbf{z}^{\mathrm{T}} \otimes \mathbf{z}^{\mathrm{T}})] - \mathrm{vec}\,\mathbf{I}_p \,\mathrm{vec}^{\mathrm{T}}\,\mathbf{I}_p)(\boldsymbol{\Sigma}^{-1/2} \otimes \boldsymbol{\Sigma}^{-1/2})G$$

$$(B.18) \quad = 1/(4q^2)c_2[(r(1-q)+1)^{-2}+1]G^{\mathrm{T}}(\boldsymbol{\Sigma}^{-1} \otimes \boldsymbol{\Sigma}^{-1})G$$

$$(B.19) \quad = 1/(4q^2)\frac{[(2(1-q)+1)^{-2}+1](3-2q)^{-p/2}}{4(2\pi q^p|\boldsymbol{\Sigma}|)^{2-q}}G^{\mathrm{T}}(\boldsymbol{\Sigma}^{-1} \otimes \boldsymbol{\Sigma}^{-1})G.$$

A similar calculation gives

$$(B.20) \quad J_{22} = 1/(4q^2)\frac{[(2-q)^{-2}+1](2-q)^{-p/2}}{(2\pi q^p|\boldsymbol{\Sigma}|)^{(2-q)/2}}G^{\mathrm{T}}(\boldsymbol{\Sigma}^{-1} \otimes \boldsymbol{\Sigma}^{-1})G.$$

Finally, we assemble (B.19) and (B.20) obtaining

$$(B.21) \quad \begin{aligned} V_{22} &= J_{22}^{-1} K_{22} J_{22}^{-1} \\ &= \frac{4q^2[(3-2q)^2+1](2-q)^{4+p}}{[(2-q)^2+1]^2(3-2q)^{2+p/2}}[G^{\mathrm{T}}(\boldsymbol{\Sigma}^{-1} \otimes \boldsymbol{\Sigma}^{-1})G]^{-1}. \end{aligned}$$

**Acknowledgments.** The authors wish to thank Tiefeng Jiang for helpful discussions. Comments from two referees, especially the one with a number of very constructive suggestions on improving the paper, are greatly appreciated.

DIPARTMENTO DI ECONOMIA POLITICA
UNIVERSITÀ DI MODENA E REGGIO EMILIA
VIA BERENGARIO 51
MODENA, 41100
ITALY
E-MAIL: davide.ferrari@unimore.it

SCHOOL OF STATISTICS
UNIVERSITY OF MINNESOTA
313 FORD HALL
224 CHURCH STREET S.E.
MINNEAPOLIS, MINNESOTA 55455
USA
E-MAIL: yyang@stat.umn.edu